\documentclass[11pt]{article}
\usepackage{amssymb,amsmath,color} % font used for R in Real numbers

 \allowdisplaybreaks

\catcode`!=11
\let\!int\int \def\int{\displaystyle\!int}
\let\!lim\lim \def\lim{\displaystyle\!lim}
\let\!sum\sum \def\sum{\displaystyle\!sum}
\let\!sup\sup \def\sup{\displaystyle\!sup}
\let\!inf\inf \def\inf{\displaystyle\!inf}
\let\!cap\cap \def\cap{\displaystyle\!cap}
\let\!max\max \def\max{\displaystyle\!max}
\let\!min\min \def\min{\displaystyle\!min}
\let\!frac\frac \def\frac{\displaystyle\!frac}
\catcode`!=12

\let\oldsection\section
\renewcommand\section{\setcounter{equation}{0}\oldsection}

\allowdisplaybreaks

\newtheorem{thm}{Theorem}[section]

\newtheorem{lem}{Lemma}[section]
\newtheorem{re}{Remark}[section]

\begin{document}
%%%%%%%%%%%%%%%%%%%% title %%%%%%%%%%%%%%%%%%%%%%%%%%%%%%%%%%%%%%%%%%%%%%%%%%
\title{Lifespan of Solution to MHD Boundary Layer Equations with Analytic Perturbation of General Shear Flow}
\author{{\bf Feng Xie}\\[1mm]
\small School of Mathematical Sciences, and LSC-MOE,\\[1mm]
\small Shanghai Jiao Tong University,
Shanghai 200240, P.R.China\\[1mm]
{\bf Tong Yang}\\[1mm]
\small Department of Mathematics,
City University of Hong Kong,\\[1mm]
\small Tat Chee Avenue, Kowloon, Hong Kong
}
\date{}
\maketitle

%\bigskip

\begin{center}
 Dedicated to Professor Philippe G. Ciarlet on the Occasion of
his 80th Birthday
\end{center}

\bigskip

\begin{abstract} In this paper, we consider the lifespan of solution to the
	MHD boundary layer system as an analytic perturbation of general shear flow.
	By using the cancellation mechanism in the system observed in
	\cite{LXY1}, the lifespan of solution is shown to have a lower bound in the order of
	$\varepsilon^{-2+}$ if the strength of the perturbation is of the order
	of $\varepsilon$. Since there is no restriction on the strength of the
	shear flow and the lifespan estimate is larger than the one obtained for the
	classical Prandtl system in this setting, it reveals the stabilizing effect
	of the magnetic field on the electrically conducting fluid near the boundary.
\end{abstract}
\footnotetext[1]{{\it E-mail address:}
tzxief@sjtu.edu.cn (F. Xie)}
\footnotetext[2]{{\it E-mail address:} matyang@cityu.edu.hk(T. Yang)}

\vskip 2mm


\vskip 2mm


\section{Introduction}
Consider the high Reynolds number limit to the MHD system near a no-slip boundary, the following MHD boundary layer system was derived in \cite{LXY1}
when both of the Reynolds number and the magnetic Reynolds number have the
same order in two space dimensions. Precisely, consider
the MHD system in  the domain $\{(x,Y)|x\in\mathbb{R}, Y\in\mathbb{R}_+\}$
with $Y=0$ being the boundary,
\begin{align}
\label{V1}
\left\{
\begin{array}{ll}
\partial_t u^\epsilon+(u^\epsilon\partial_x+v^\epsilon\partial_Y)u^\epsilon+\partial_xp^\epsilon-(h^\epsilon\partial_x+g^\epsilon\partial_Y)h^\epsilon=\epsilon(\partial_x^2 u^\epsilon+\partial_Y^2 u^\epsilon),\\
\partial_t v^\epsilon+(u^\epsilon\partial_x+v^\epsilon\partial_Y)v^\epsilon+\partial_Yp^\epsilon-(h^\epsilon\partial_x+g^\epsilon\partial_Y)g^\epsilon=\epsilon(\partial_x^2 v^\epsilon+\partial_Y^2 v^\epsilon),\\
\partial_t h^\epsilon+(u^\epsilon\partial_x+v^\epsilon\partial_Y)h^\epsilon-(h^\epsilon\partial_x+g^\epsilon\partial_Y)u^\epsilon=\kappa\epsilon(\partial_x^2 h^\epsilon+\partial_Y^2 h^\epsilon),\\
\partial_t g^\epsilon+(u^\epsilon\partial_x+v^\epsilon\partial_Y)g^\epsilon-(h^\epsilon\partial_x+g^\epsilon\partial_Y)v^\epsilon=\kappa\epsilon(\partial_x^2 g^\epsilon+\partial_Y^2 g^\epsilon),\\
\partial_x u^\epsilon+\partial_Y v^\epsilon=0,\qquad \partial_x h^\epsilon+\partial_Y g^\epsilon=0,
\end{array}
\right.
\end{align}
where both the viscosity and resistivity coefficients are denoted
by a small positive parameter $\epsilon$, $(u^\epsilon, v^\epsilon)$
and $(h^\epsilon, g^\epsilon)$ represent the velocity and
 the magnetic field respectively.
 The no-slip boundary condition is imposed on the velocity field
\begin{align}
\label{BCV}
(u^\epsilon, v^\epsilon)|_{Y=0}={\bf{0}},
\end{align}
and the perfectly conducting boundary condition is given for the magnetic
 field
\begin{align}
\label{BCM}
(\partial_Yh^\epsilon, g^\epsilon)|_{Y=0}={\bf{0}}.
\end{align}
Formally, when $\epsilon=0$,  (\ref{V1}) is reduced into the following incompressible ideal MHD system
\begin{align}
\label{ILE0}
\left\{
\begin{array}{ll}
\partial_t u^0_e+(u^0_e\partial_x+v^0_e\partial_Y)u^0_e+\partial_xp^0_e-(h^0_e\partial_x+g^0_e\partial_Y)h^0_e=0,\\
\partial_t v^0_e+(u^0_e\partial_x+v^0_e\partial_Y)v^0_e+\partial_Yp^0_e-(h^0_e\partial_x+g^0_e\partial_Y)g^0_e=0,\\
\partial_t h^0_e+(u^0_e\partial_x+v^0_e\partial_Y)h^0_e-(h^0_e\partial_x+g^0_e\partial_Y)u^0_e=0,\\
\partial_t g^0_e+(u^0_e\partial_x+v^0_e\partial_Y)g^0_e-(h^0_e\partial_x+g^0_e\partial_Y)v^0_e=0,\\
\partial_x u^0_e+\partial_Y v^0_e=0,\qquad \partial_x h^0_e+\partial_Y g^0_e=0.
\end{array}
\right.
\end{align}
Since the sovability of the system (\ref{ILE0}) requires only the
 normal components of the velocity and magnetic fields $(v_e^0, g^0_e)$ on the boundary
\begin{align}
\label{IBE}
(v^0_e, g^0_e)|_{Y=0}={\bf 0},
\end{align}
in the limit from (\ref{V1}) to (\ref{ILE0}), a Prandtl-type boundary layer can be derived to resolve the mis-match of the tangential components between the viscous flow $(u^\epsilon, h^\epsilon)$ and invicid flow $(u^0, h^0)$ on the  boundary $\{Y=0\}$. And this system governing the fluid behavior in the
leading order of approximation near
the boundary is  derived in \cite{GP,LXY1,LXY2}:
\begin{align}
\label{1.1}
\left\{
\begin{array}{ll}
\partial_tu_1+u_1\partial_xu_1+u_2\partial_yu_1=b_1\partial_xb_1+b_2\partial_yb_1+\partial^2_yu_1,\\
\partial_tb_1+\partial_y(u_2b_1-u_1b_2)=\kappa\partial_y^2b_1,\\
\partial_xu_1+\partial_yu_2=0,\quad \partial_xb_1+\partial_yb_2=0
\end{array}
\right.
\end{align}
in $\mathbb{H}=\{(x,y)\in\mathbb{R}^2| y\geq0\}$ with the fast variable $y=Y/\sqrt{\epsilon}$. Here, the trace of the horizontal ideal MHD flow (\ref{ILE0}) on the boundary $\{Y=0\}$ is
 assumed to be a constant vector so that the pressure term $\partial_xp^0_e(t,x,0)$ vanishes by the Bernoulli's law.

Consider the system $(\ref{1.1})$ with initial data
\begin{align}
\label{ID}
u_1(t,x,y)|_{t=0}=u_0(x,y),\qquad\qquad b_1(t,x,y)|_{t=0}=b_0(x,y),
\end{align}
and the boundary conditions
\begin{align}
\label{1.2}
\left\{
\begin{array}{ll}
u_1|_{y=0}=0,\\
u_2|_{y=0}=0,
\end{array}
\right.
\qquad\hbox{and}\qquad
\left\{
\begin{array}{ll}
\partial_yb_1|_{y=0}=0,\\
b_2|_{y=0}=0.
\end{array}
\right.
\end{align}
%These boundary conditions correspond to the no-slip boundary condition of velocity (\ref{BCV}) and perfectly conducting boundary condition of magnetic fields (\ref{BCM}) for 2D incompressible MHD equations (\ref{V1}).
And the far field state is denoted by $(\bar{u}, \bar{b})$:
\begin{align}
\label{1.3}
\lim_{y\rightarrow+\infty}u_1=u^0_e(t,x,0)\triangleq\bar{u},\qquad \lim_{y\rightarrow+\infty}b_1=h^0_e(t,x,0)\triangleq\bar{b}.
\end{align}
First of all, a  shear flow $(u_s(t,y),0,\bar{b}, 0)$ is a trivial
solution to the system  \eqref{1.1} with $u_s(t,y)$ solving
\begin{align}
\label{GS}
\left\{
\begin{array}{ll}
\partial_tu_s(t,y)-\partial_y^2u_s(t,y)=0, \quad (t, y)\in\mathbb{R}_+\times\mathbb{R}_+,\\
u_s(t,y=0)=0, \qquad \lim_{y\rightarrow\infty}u_s(t,y)=\bar{u},\\
u_{s}(t=0, y)=u_{s0}(y).
\end{array}
\right.
\end{align}
In the following discussion, we assume the shear flow $u_s(t,y)$ has the following properties:
\begin{align*}
{\bf (H)}\quad \|\partial_y^iu_s(t,\cdot)\|_{L^\infty_y}\leq \frac{C}{\langle t \rangle^{i/2}}\quad (i=1,2),\quad \int_0^\infty|\partial_yu_s(t,y)|dy<C,\quad\|\theta_\alpha\partial_y^2u_s(t,\cdot)\|_{L^2_y}\leq \frac{C}{\langle t \rangle^{3/4}},
\end{align*}
for some generic constant $C$.
\begin{re}
	\label{RK2}
	The assumption (H) on the shear flow holds
	for  a large class of initial data $u_{s0}$.
	For example, it holds for the initial data $u_{s0}=\chi(y)$ with $\chi(y)\in C^\infty(\mathbb{R})$, $\chi(y)=0$ for $y\leq 1$ and $\chi(y)=\bar{u}$ for $y\geq 2$ considered in \cite{ZZ} for the Prandtl system. Note that here we do
	not assume the smallness of the shear flow. In addition, it also holds when  $u_{s0}(y)=\frac{1}{\sqrt{\pi}}\int_0^{y}\exp(-\frac{z^2}{4})dz$
	 considered in \cite{IV} for the Prandtl system where the almost global solution is
	 obtained. Note that for the classical Prandtl equations, such shear flow
	 in the form of Guassian error
	 function yields a time decay damping term in the time evolution equation of
	 $u_1$, however, it does not leads to  any damping effect in the MHD boundary
	 layer system \eqref{1.1}.
\end{re}

To define the function space of the solution considered in this paper,
the following Gaussian weighted function $\theta_\alpha$ will be used:
\begin{align*}
\theta_{\alpha}(t,y)=\exp{(\frac{\alpha z(t,y)^2}{4})},\quad \hbox{with}\quad z(t,y)=\frac{y}{\sqrt{\langle t\rangle}},\quad \langle t\rangle=1+t\quad\hbox{and}\quad \alpha\in[1/4,1/2].
\end{align*}
With this and
\begin{align*}
M_m=\frac{\sqrt{m+1}}{m!},
\end{align*}
define  the Sobolev weighted semi-norms by
\begin{align}
\label{NM1}
&X_m=X_m(f,\tau)=\|\theta_{\alpha}\partial_x^mf\|_{L^2}\tau^mM_m,\quad
D_m=D_m(f,\tau)=\|\theta_{\alpha}\partial_y\partial_x^mf\|_{L^2}\tau^mM_m,\nonumber\\
&Z_m=Z_m(f,\tau)=\|z\theta_{\alpha}\partial_x^mf\|_{L^2}\tau^mM_m,\quad
Y_m=Y_m(f,\tau)=\|\theta_{\alpha}\partial_x^mf\|_{L^2}\tau^{m-1}mM_m.
\end{align}
Then the following space of analytic functions
in the tangential variable $x$ and Sobolev weighted in the normal variable $y$
is defined by
\begin{align*}
X_{\tau, \alpha}=\{f(t,x,y)\in L^2(\mathbb{H}; \theta_{\alpha}dxdy): \|f\|_{X_{\tau, \alpha}}<\infty\}
\end{align*}
with $\tau>0$ and the norm
\begin{align*}
\|f\|_{X_{\tau, \alpha}}=\sum_{m\geq 0}X_m(f,\tau).
\end{align*}
In addition,  the following two semi-norms will also be used:
\begin{align*}
\|f\|_{D_{\tau, \alpha}}=\sum_{m\geq 0}D_m(f,\tau)=\|\partial_yf\|_{X_{\tau, \alpha}},\quad
\|f\|_{Y_{\tau, \alpha}}=\sum_{m\geq 1}Y_m(f,\tau).
\end{align*}
Here, the summation over $m$ is considered in the $l^1$ sense that is similar
to the definition used in \cite{IV,ZZ} rather than in the $l^2$ sense used in \cite{KV}. With the above notations, we are now
 ready to state the main Theorem as follows.

\begin{thm}
\label{THM}
For any $\lambda\in [3/2, 2)$, there exists a small positive constant
$\varepsilon_*$ depending on $2-\lambda$.
Under the assumption (H) on the backgroud shear flow
$(u_s(t,y), 0, \bar{b},0)$ with $\bar{b}\neq 0$, assume the initial data $u_0$ and $b_0$  satisfy
\begin{align}
\label{THM1}
\|u_0-u_s(0,y)\|_{X_{2\tau_0,1/2}}\leq \varepsilon,
\quad
\|b_0-\bar{b}\|_{X_{2\tau_0,1/2}}\leq \varepsilon,
\end{align}
for some given $\varepsilon\in (0, \varepsilon_*]$. Then there exists a unique solution $(u_1, u_2, b_1, b_2)$ to the MHD boundary layer equations (\ref{1.1})-(\ref{1.3}) such that
\begin{align*}
(u_1-u_s(t,y),  b_1-\bar{b})\in X_{\tau,\alpha},\  \alpha\in[1/4,1/2],
\end{align*}
with analyticity radius $\tau$ larger than $\tau_0/4$ in the time interval $[0, T_\varepsilon]$.
And the lifespan $T_\varepsilon$ has the following low bound estimate
\begin{align}
\label{THM3}
T_\varepsilon\geq {C}\varepsilon^{-\lambda},
\end{align}
where the constant ${C}$ is  independent of $\varepsilon$.
%It is also emphasized that $\varepsilon$ and $\tau_0$ are independent of $\bar{u}$ and $\bar{b}$.
\end{thm}
%\begin{re}
%\label{RK1}
%The definitions of function spaces $X_{\tau,\alpha},  \alpha\in[1/4,1/2]$  will be given in details in Section 2. Where we also introduce other analytic function spaces used in this paper.

%It is noted that the index $3/2$ is not essential. When the index  is fixed, we can determine the sufficiently small universal constant $\varepsilon_*>0$ in an explicitly form. Moreover, when the index $\lambda$ is more close to $2$, then we require the universal constant $\varepsilon_*$ is more small accordingly.
%\end{re}

As is well-known that the leading order characteristic boundary layer for the
incompressible Navier-Stokes equations
with no-slip boundary condition is described by the classical Prandtl
equations derived by Prandtl  \cite{P} in 1904. In the two space dimensions,
under the monotonicity
assumption on the tangential velocity in the normal direction, Oleinik firstly
obtained the local existence of classical solutions
by using the Crocco transformation, cf.  \cite{O} and
Oleinik-Samokhin's classical
book \cite{OS}. Recently, this well-posedness result was re-proved by using an energy method  in the framework of
 Sobolev spaces in \cite{AWXY} and \cite{MW1} independently by observing the cancellation mechanism in the convection terms. And by imposing an additional
favorable condition on the pressure,
a global in time weak solution was obtained in
\cite{XZ}.

When the monotonicity condition is violated, singularity formation or separation of the boundary
layer is well expected and observed. For this, E-Engquist constructed a finite
time blowup solution to the Prandtl
equations in \cite{EE}. Recently, when the background shear flow has a non-degenerate critical point, some interesting ill-posedness (or instability) phenomena of solutions to both linear and nonlinear classical Prandtl equations around
shear flows are studied, cf.
\cite{GD,GN,GN1}. All these results show that the monotonicity assumption on the tangential velocity plays a key role for well-posedness theory
except in the frameworks of analytic functions and Gevrey regularity classes. Indeed, in the framework of
analytic functions, Sammartino and Caflisch \cite{SC,CS} established
the local well-posedness theory of the Prandtl system in three space dimensions and also justified the Prandtl ansatz in this setting by applying
the abstract Cauchy-Kowalewskaya (CK) theorem initated by Asano's unpublished work.
 Later, the analyticity requirement in the normal variable $y$ was removed by
Lombardo, Cannone and Sammartino in \cite{LCS} because of the viscous effect in the normal direction.

Recently, Zhang and Zhang obtained the lifespan of small analytic solution to the classical Prandtl equations with small analytic initial data in \cite{ZZ}. Precisely, when the strength of background shear flow is of the order
 of $\varepsilon^{5/3}$ and the perturbation is of the order of $\varepsilon$, they showed that  the classical Prandtl
system has a unique solution with a lower bound estimate on the lifespan in the order of  $\varepsilon^{-4/3}$.
Furthermore, if the initial data is a small analytic perturbation of the Guassian error function (\ref{GS}),  an
almost global existence for the Prandtl boundary layer equations is proved in \cite{IV}.

On the other hand, to study the high Reynolds number limits for  the MHD equations (\ref{V1}) with no-slip boundary condition on the velocity (\ref{BCV}) and perfect conducting boundary condition (\ref{BCM}) on the magnetic field,
one can apply the Prandtl ansatz to derive the boundary layer system  (\ref{1.1}) as
 the leading order description on the flow near the boundary. For this, readers can
  refer to \cite{GP,LXY1,LXY2,LXY3,XY} about  the formal derivation of (\ref{1.1}), the well-posedness theory
  of the system  and the justification of the Prandtl ansatz locally in time.

This paper is about  long time existence of solutions to (\ref{1.1})-(\ref{1.3}). Precisely, we will show that if the initial data is a small perturbation of a shear flow analytically in the order of  $\varepsilon$, then there exists a unique solution to (\ref{1.1})-(\ref{1.3}) with the lifespan $T_\varepsilon$ of the order of $ \varepsilon^{-2+}$.  Compared with the
estimate on the lifespan of solutions to the classical Prandtl system studied in  \cite{ZZ}, the lower bound estimate
is larger and there is no requirement on the smallness of the background shear flow because the mechanism in the
system is used due to the non-degeneracy of the tangential  magnetic field. However, it is not known
whether one can obtain a global or
almost global in time solution like the work on the Prandtl system  when the background shear velocity is taken to be a Gaussian error function in \cite{IV}. We mention that even though  Lin and Zhang showed the almost global existence of solution to MHD boundary layer equations with zero Dirichlet boundary condition on the  magnetic field in \cite{LZ} when
the components of both the background velocity and magnetic fields are Guassian error functions, it is not clear wheather
 the system (\ref{1.1}) holds with zero Dirichlet boundary condition even in formal derivation.

The analysis on the lifespan of the perturbed system in this paper relies on the introduction of some new unknown functions that capture the cancellation of some  linear terms. Unlike the work
in \cite{IV} on the  Prandtl system for which the cancellation
yields  a damping  term in the time evolution of the perturbation of the tangetial velocity field, there
is no such damping effect observed for the  MHD boundary layer system.

Finally, the rest of the paper is organized as follows. After giving  some preliminary estimates,  a
uniform estimate on the solution will be proved in the next section.  Based on this uniform estimate,  a low bound of the lifespan of  solution is  derived  in Section 3. The uniqueness part is done in Section 4. Throughout the paper,  constants denoted by  $C$, $\bar{C}, C_0, C_1$ and $C_2$ are generic and independent of the small parameter $\varepsilon$.

\section{Uniform Estimate}

We first list the following two priliminary estimates on the functions in the norms defined in the previous section. The first
estimate indeed is from  Lemma 3.3 in \cite{IV} (also see \cite{H}).
\begin{lem}
\label{LEM2.1}
(Poincar\'e type inequality with Gaussian weight) Let $f$ be a function such that $f|_{y=0}=0\ (or\ \partial_yf|_{y=0})$ and $f|_{y=\infty}=0$. Then, for $\alpha\in[1/4,1/2], m\geq 0$ and $t\geq 0$, it holds that
\begin{align}
\label{2.1}
\frac{\alpha}{\langle t\rangle}\|\theta_{\alpha}\partial_x^mf\|_{L^2_y}^2\leq\|\theta_{\alpha}\partial_y\partial_x^mf\|_{L^2_y}^2.
\end{align}
\end{lem}

The second lemma is  used  in \cite{IV} and we include it here with a short proof for convenience of readers.

\begin{lem}\label{LEM2.2}
Let $f$ be a function such that $f|_{y=0}=0\ (or\ \partial_yf|_{y=0})$ and $f|_{y=\infty}=0$.  Then
\begin{align}
\label{LL}
\sum_{m\geq 0}\frac{\|\theta_{\alpha}\partial_y\partial_x^mf\|^2_{L^2}}{\|\theta_{\alpha}\partial_x^mf\|_{L^2}}\tau^mM_m\geq \frac{\alpha^{1/2}\beta}{2\langle t\rangle^{1/2}}\|f\|_{D_{\tau,\alpha}}+\frac{\alpha(1-\beta)}{\langle t\rangle}\|f\|_{X_{\tau,\alpha}},
\end{align}
for $\beta\in(0, 1/2)$.
\end{lem}
\begin{pf}
In fact, by Lemma \ref{LEM2.1}, one has
\begin{align*}
\frac{\|\theta_{\alpha}\partial_y\partial_x^mf\|^2_{L^2}}{\|\theta_{\alpha}\partial_x^mf\|_{L^2}}\geq & \frac{\beta}{2}\frac{\|\theta_{\alpha}\partial_y\partial_x^mf\|^2_{L^2}}{\|\theta_{\alpha}\partial_x^mf\|_{L^2}}+\frac{2-\beta}{2}\frac{\alpha^{1/2}}{\langle t\rangle^{1/2}}\|\theta_{\alpha}\partial_y\partial_x^mf\|_{L^2}\\
\geq &\frac{\beta}{2}\frac{\|\theta_{\alpha}\partial_y\partial_x^mf\|^2_{L^2}}{\|\theta_{\alpha}\partial_x^mf\|_{L^2}}+\frac{\beta\alpha^{1/2}}{2\langle t\rangle^{1/2}}\|\theta_{\alpha}\partial_y\partial_x^mf\|_{L^2}+\frac{\alpha(1-\beta)}{\langle t\rangle}\|\theta_{\alpha}\partial_x^mf\|_{L^2}\\
\geq& \frac{\beta\alpha^{1/2}}{2\langle t\rangle^{1/2}}\|\theta_{\alpha}\partial_y\partial_x^mf\|_{L^2}+\frac{\alpha(1-\beta)}{\langle t\rangle}\|\theta_{\alpha}\partial_x^mf\|_{L^2}.
\end{align*}
Multiplying the above inequality by $\tau^mM_m$  and summing up in $m\geq 0$ give (\ref{LL}).
\end{pf}

 %\section{Uniform {\it A Priori} Estimates}
 We are now ready to study a uniform estimate on the solution. For this, we first
rewrite the solution to (\ref{1.1})-(\ref{1.3}) as a perturbation $(u, v, b, g)$ of the $(u_s(t,y), 0, \bar{b}, 0)$ by denoting
\begin{align}
\label{EXP}
\left\{
\begin{array}{ll}
u_1=u_s(t,y)+u,\\
u_2=v,
\end{array}
\right.
\qquad\qquad
\left\{
\begin{array}{ll}
b_1=\bar{b}+b,\\
b_2=g.
\end{array}
\right.
\end{align}
Without loss of generality, take $\bar{b}=1$
and $\kappa=1$. Then  (\ref{1.1}) yields
\begin{align}
\label{1.4}
\left\{
\begin{array}{ll}
\partial_tu+(u_s+u)\partial_xu+v\partial_y(u_s+u)-(1+b)\partial_xb-g\partial_yb-\partial_y^2u=0,\\
\partial_tb-(1+b)\partial_xu-g\partial_y(u_s+u)+(u_s+u)\partial_xb+v\partial_yb-\partial_y^2b=0.
\end{array}
\right.
\end{align}
And the initial and boundary data of $(u, v)$ and $(b, g)$ are given by
\begin{align}
\label{NID}
u(t,x,y)|_{t=0}=u_0(x,y)-u_s(0,y),\qquad b(t, x,y)|_{t=0}=b_0(x,y)-1,
\end{align}
\begin{align}
\label{NBC}
\left\{
\begin{array}{ll}
u|_{y=0}=0,\\
v|_{y=0}=0,
\end{array}
\right.
\qquad\hbox{and}\qquad
\left\{
\begin{array}{ll}
\partial_yb|_{y=0}=0,\\
g|_{y=0}=0,
\end{array}
\right.
\end{align}
with the corresponding far field condition
\begin{align}
\label{NFC}
\lim_{y\rightarrow+\infty}u=0,\qquad \lim_{y\rightarrow+\infty}b=0.
\end{align}

It suffices to establish the long time existence of solutions to (\ref{1.4})-(\ref{NFC}). In this section, we focus on the uniform a priori estimate on  the solution to (\ref{1.4}) in the analytical framework defined in Section 1.

Integrating equation $(\ref{1.4})_2$ over $[0,y]$ gives that
\begin{align}
\label{3.1}
\partial_t\int_0^ybd\tilde{y}+v(1+b)-(u_s+u)g=\partial^2_y\int_0^ybd\tilde{y},
\end{align}
where the boundary conditions that $\partial_yb|_{y=0}=v|_{y=0}=g|_{y=0}=0$ are used.\\
Define
\begin{align*}
\psi(t,y)=\int_0^ybd\tilde{y},
\end{align*}
one has
\begin{align}
\label{3.2}
\partial_t\psi+v(1+b)-(u_s+u)g=\partial^2_y\psi.
\end{align}
Now introduce  new unknown functions by taking care of the cancellation mechamism in the system as obseved in \cite{LXY1} as follows
\begin{align}
\label{3.3}
\tilde{u}=u-\partial_yu_s\psi,\qquad \tilde{b}=b.
\end{align}
Then $(\tilde{u}, \tilde{b})$ satisfies the following equations.
\begin{align}
\label{3.4}
\left\{
\begin{array}{ll}
\partial_t\tilde{u}-\partial_y^2\tilde{u}+(u_s+u)\partial_x\tilde{u}+v\partial_y\tilde{u}-(1+b)\partial_x\tilde{b}-g\partial_y\tilde{b}-2\partial_y^2u_s\tilde{b}+v\partial_y^2u_s\psi=0,\\
\partial_t\tilde{b}-\partial_y^2\tilde{b}-(1+b)\partial_x\tilde{u}-g\partial_y\tilde{u}+(u_s+u)\partial_x\tilde{b}+v\partial_y\tilde{b}-g\partial_y^2u_s\psi=0.
\end{array}
\right.
\end{align}
Here we have used the following fact that $u_s$ is the solution to the heat equation. That is,
\begin{align*}
\partial_tu_s-\partial_y^2u_s=0,\qquad \partial_t\partial_yu_s-\partial_y^3u_s=0.
\end{align*}
By a direct calculation, the boundary conditions of $(\tilde{u}, \tilde{b})$ are given by
\begin{align}
\label{BC}
\tilde{u}|_{y=0}=0,\qquad \partial_y\tilde{b}|_{y=0}=0,
\end{align}
\begin{align}
\label{FBC}
\tilde{u}|_{y=\infty}=0,\qquad \tilde{b}|_{y=\infty}=0.
\end{align}
We then turn to show the existence of solution $(\tilde{u}, \tilde{b})$ to (\ref{3.4})-(\ref{FBC}) with the
corresponding initial data.
\begin{align}
\label{IIII}
\tilde{u}(0,x,y)=u(0,x,y)-\partial_yu_s(0,y)\int_0^yb(0,x,\tilde{y})d\tilde{y},\qquad \tilde{b}(0,x,y)=b(0,x,y).
\end{align}
Note that
\begin{align}
\label{IIIIE}
\|\tilde{u}(0,x,y)\|_{X_{2\tau_0, \alpha}}\leq \|u(0,x,y)\|_{X_{2\tau_0, \alpha}}+C\|b(0,x,y)\|_{X_{2\tau_0, \alpha}},
\end{align}
for $\alpha\in[1/4, 1/2]$.

Moreover, once the existence of solution $(\tilde{u}, \tilde{b})$ to (\ref{3.4})-(\ref{IIII}) is obtained,  one can define $(u, b)$ by
\begin{align}
\label{RS}
u(t,x,y)=\tilde{u}(t,x,y)+\partial_yu_s(t,y)\int_0^y\tilde{b}(t,x,\tilde{y})d\tilde{y},\qquad b(t,x,y)=\tilde{b}(t,x,y).
\end{align}
It is straightforward to check that $(u, b)$ is a solution to (\ref{1.4})-(\ref{NFC}) with the following estimates
\begin{align*}
\|u\|_{X_{\tau, \alpha}}\leq \|\tilde{u}\|_{X_{\tau, \alpha}}+C\|\tilde{b}\|_{X_{\tau, \alpha}},\qquad  \|b\|_{X_{\tau, \alpha}}=\|\tilde{b}\|_{X_{\tau, \alpha}}.
\end{align*}
Therefore, we only need to estimate  the solution $(\tilde{u}, \tilde{b})$ to (\ref{3.4})-(\ref{IIII})
in the analytic norms as shown in the next two subsections.

\subsection{A priori estimate on velocity field}
For  $m\geq 0$,  by applying the tangential derivative operator $\partial_x^m$ to $(\ref{3.4})_1$ and multiplying it  by $\theta_{\alpha}^2\partial_x^m\tilde{u}$, the integration over  $\mathbb{H}$ yields
\begin{align}
\label{3.5}
\int_{\mathbb{H}}\partial_x^m(\partial_t\tilde{u}-\partial_y^2\tilde{u}+(u_s+u)\partial_x\tilde{u}+v\partial_y\tilde{u}-
(1+b)\partial_x\tilde{b}-g\partial_y\tilde{b}-2\partial_y^2u_s\tilde{b}+v\psi\partial_y^2u_s)\theta_{\alpha}^2\partial_x^m\tilde{u}dxdy=0.
\end{align}
We now  estimate each term in $(\ref{3.5})$ as follows. Firstly, note that
\begin{align}
\label{3.6}
&\int_\mathbb{H}\partial_t\partial_x^m\tilde{u}\theta_{\alpha}^2\partial_x^m\tilde{u}dxdy\nonumber\\
= &\frac12\frac{d}{dt}\int_\mathbb{H}(\partial_x^m\tilde{u})^2\theta_{\alpha}^2dxdy-\int_\mathbb{H}(\partial_x^m\tilde{u})^2\theta_{\alpha}\frac{d}{dt}\theta_{\alpha}dxdy\\
= &\frac12\frac{d}{dt}\|\theta_{\alpha}\partial_x^m\tilde{u}\|^2_{L^2}+\frac{\alpha}{4\langle t\rangle}\|\theta_{\alpha}z\partial_x^m\tilde{u}\|^2_{L^2},\nonumber
\end{align}
and
\begin{align*}
-\int_\mathbb{H}\partial_y^2\partial_x^m\tilde{u}\theta_{\alpha}^2\partial_x^m\tilde{u}dxdy
= \|\theta_{\alpha}\partial_x^m\partial_y\tilde{u}\|_{L^2}^2+\int_\mathbb{H}\partial_y\partial_x^m\tilde{u}\partial_y(\theta_{\alpha}^2)\partial_x^m\tilde{u}dxdy.
\end{align*}
The boundary  term vanishes because  of the boundary condition  $\partial_x^m\tilde{u}|_{y=0}=0$. Furthermore,
\begin{align*}
&\int_\mathbb{H}\partial_y\partial_x^m\tilde{u}\partial_y(\theta_{\alpha}^2)\partial_x^m\tilde{u}dxdy
=-\frac12\int_\mathbb{H}(\partial_x^m\tilde{u})^2\partial_y^2(\theta_{\alpha}^2)dxdy\\
=&-\frac{\alpha}{2}\frac{1}{\langle t\rangle}\|\theta_{\alpha}\partial_x^m\tilde{u}\|^2_{L^2}-\frac{\alpha^2}{2}\frac{1}{\langle t\rangle}\|\theta_{\alpha}z\partial_x^m\tilde{u}\|^2_{L^2},
\end{align*}
where we have used
\begin{align*}
\partial_y^2(\theta_{\alpha}^2)=\frac{\alpha}{\langle t\rangle}\theta_{\alpha}^2+\frac{\alpha^2}{\langle t\rangle}z^2(t,y)\theta_{\alpha}^2.
\end{align*}
Consequently,
\begin{align}
\label{3.7}
-\int_\mathbb{H}\partial_y^2\partial_x^m\tilde{u}\theta_{\alpha}^2\partial_x^m\tilde{u}dxdy
= \|\theta_{\alpha}\partial_x^m\partial_y\tilde{u}\|_{L^2}^2-\frac{\alpha}{2}\frac{1}{\langle t\rangle}\|\theta_{\alpha}\partial_x^m\tilde{u}\|^2_{L^2}-\frac{\alpha^2}{2}\frac{1}{\langle t\rangle}\|\theta_{\alpha}z\partial_x^m\tilde{u}\|^2_{L^2}.
\end{align}
For the nonlinear terms in (\ref{3.5}), we have
\begin{align*}
\int_{\mathbb{H}}\partial_x^m((u_s+u)\partial_x\tilde{u})\theta_{\alpha}^2\partial_x^m\tilde{u}dxdy
=\sum_{j=0}^m(\begin{array}{ll}m\\j\end{array})\int_{\mathbb{H}}\partial^{m-j}_xu\partial_x^{j+1}\tilde{u}\theta_{\alpha}^2\partial_x^m\tilde{u}dxdy
\triangleq R_1
\end{align*}
and
\begin{align*}
|R_1|\leq& \sum_{j=0}^{[m/2]}(\begin{array}{ll}m\\j\end{array})\|\partial^{m-j}_xu\|_{L^2_xL^\infty_y}\|\theta_{\alpha}\partial^{j+1}_x\tilde{u}\|_{L^\infty_xL^2_y}\|\theta_{\alpha}\partial^{m}_x\tilde{u}\|_{L^2}\\
+&\sum_{j=[m/2]+1}^{m}(\begin{array}{ll}m\\j\end{array})\|\partial^{m-j}_xu\|_{L^\infty_{xy}}\|\theta_{\alpha}\partial^{j+1}_x\tilde{u}\|_{L^2}\|\theta_{\alpha}\partial^{m}_x\tilde{u}\|_{L^2}.
\end{align*}
For $0\leq j\leq [m/2]$, by (\ref{RS}), one has
\begin{align*}
&\|\partial^{m-j}_xu\|_{L^2_xL^\infty_y}=\|\partial^{m-j}_x(\tilde{u}+\partial_yu_s\psi)\|_{L^2_xL^\infty_y}\\
\leq &\|\partial^{m-j}_x\tilde{u}\|_{L^2_xL^\infty_y}+\|\partial_yu_s\partial^{m-j}_x\psi\|_{L^2_xL^\infty_y}\\
\leq &C\|\theta_{\alpha}\partial^{m-j}_x\tilde{u}\|^{1/2}_{L^2}\|\theta_{\alpha}\partial^{m-j}_x\partial_y\tilde{u}\|^{1/2}_{L^2}+C\langle t\rangle^{-1/4}\|\theta_{\alpha}\partial^{m-j}_x\tilde{b}\|_{L^2},
\end{align*}
where in the last inequality, we have used
$$
\|\partial_yu_s\|_{L_y^\infty}\leq \frac{C}{\sqrt{\langle t\rangle}},
$$
according to the assumption (H). Moreover,
\begin{align*}
&\|\partial^{m-j}_x\psi\|_{L^2_xL^\infty_y}=\|\int_0^y\partial^{m-j}_x\tilde{b}d\tilde{y}\|_{L^2_xL^\infty_y}\\
=&\|\int_0^y\theta_{\alpha}\partial^{m-j}_x\tilde{b}\exp(-\frac{\alpha}{4}z^2)d\tilde{y}\|_{L^2_xL^\infty_y}
\leq C\langle t\rangle^{1/4}\|\theta_{\alpha}\partial^{m-j}_x\tilde{b}\|_{L^2}.
\end{align*}
And
\begin{align*}
\|\theta_{\alpha}\partial^{j+1}_x\tilde{u}\|_{L^\infty_xL^2_y}\leq C\|\theta_{\alpha}\partial^{j+1}_x\tilde{u}\|^{1/2}_{L^2}\|\theta_{\alpha}\partial^{j+2}_x\tilde{u}\|^{1/2}_{L^2}.
\end{align*}
For $[m/2]+1\leq j\leq m$, we have
\begin{align*}
\|\partial^{m-j}_xu\|_{L^\infty_{xy}}\leq &\|\partial^{m-j}_x\tilde{u}\|_{L_{xy}^\infty}+\|\partial_yu_s\partial^{m-j}_x\psi\|_{L_{xy}^\infty}\\
\leq &C\|\theta_{\alpha}\partial^{m-j}_x\tilde{u}\|^{1/4}_{L^2}\|\theta_{\alpha}\partial^{m-j}_x\partial_y\tilde{u}\|^{1/4}_{L^2}
\|\theta_{\alpha}\partial^{m-j+1}_x\tilde{u}\|^{1/4}_{L^2}\|\theta_{\alpha}\partial^{m-j+1}_x\partial_y\tilde{u}\|^{1/4}_{L^2}\\
&+C\langle t\rangle^{-1/4}\|\theta_{\alpha}\partial^{m-j}_x\tilde{b}\|^{1/2}_{L^2}\|\theta_{\alpha}\partial^{m-j+1}_x\tilde{b}\|^{1/2}_{L^2}.
\end{align*}
Hence,
\begin{align}
\label{NE1}
\frac{|R_1|\tau^mM_m}{\|\theta_{\alpha}\partial_x^m\tilde{u}\|_{L^2}}\leq& \frac{C}{(\tau(t))^{1/2}}\left\{\sum_{j=0}^{[m/2]}(X^{1/2}_{m-j}D^{1/2}_{m-j}+\langle t\rangle^{-1/4}\bar{X}_{m-j})Y^{1/2}_{j+1}Y^{1/2}_{j+2} \right.\\
&\left.+\sum_{j=[m/2]+1}^{m}(X^{1/4}_{m-j}X^{1/4}_{m-j+1}D^{1/4}_{m-j}D^{1/4}_{m-j+1}+\langle t\rangle^{-1/4}\bar{X}_{m-j}^{1/2}\bar{X}_{m-j+1}^{1/2})Y_{j+1}\right\}.\nonumber
\end{align}
From now on, we use $X_i, D_i, Y_i$ to denote the semi-norms of function $\tilde{u}$ defined in (\ref{NM1}), and
 $\bar{X}_i, \bar{D}_i$ and $\bar{Y}_i$ for the corresponding semi-norms for $\tilde{b}$. Note that
\begin{align*}
\int_{\mathbb{H}}\partial_x^m(v\partial_y\tilde{u})\theta_{\alpha}^2\partial_x^m\tilde{u}dxdy
=\sum_{j=0}^m(\begin{array}{ll}m\\j\end{array})\int_{\mathbb{H}}\partial^{m-j}_xv\partial_x^{j}\partial_y\tilde{u}\theta_{\alpha}^2\partial_x^m\tilde{u}dxdy
\triangleq R_2
\end{align*}
and
\begin{align*}
|R_2|\leq& \sum_{j=0}^{[m/2]}(\begin{array}{ll}m\\j\end{array})\|\partial^{m-j}_xv\|_{L^2_xL^\infty_y}\|\theta_{\alpha}\partial^{j}_x\partial_y\tilde{u}\|_{L^\infty_xL^2_y}\|\theta_{\alpha}\partial^{m}_x\tilde{u}\|_{L^2}\\
+&\sum_{j=[m/2]+1}^{m}(\begin{array}{ll}m\\j\end{array})\|\partial^{m-j}_xv\|_{L^\infty_{xy}}\|\theta_{\alpha}\partial^{j}_x\partial_y\tilde{u}\|_{L^2}\|\theta_{\alpha}\partial^{m}_x\tilde{u}\|_{L^2}.
\end{align*}
For $0\leq j\leq [m/2]$,
\begin{align*}
\|\partial^{m-j}_xv\|_{L^2_xL^\infty_y}=&\|\int_0^y\partial^{m-j+1}_xud\tilde{y}\|_{L^2_xL^\infty_y}\\
\leq&\|\int_0^y\partial^{m-j+1}_x\tilde{u}d\tilde{y}\|_{L^2_xL^\infty_y}+\|\int_0^y\partial_yu_s(\int_0^{\tilde{y}}\partial^{m-j+1}_x\tilde{b}ds)d\tilde{y}\|_{L^2_xL^\infty_y}\\
\leq&\|\int_0^y\partial^{m-j+1}_x\tilde{u}d\tilde{y}\|_{L^2_xL^\infty_y}+\|\int_0^y\partial_yu_s d\tilde{y}\|_{L^\infty_y}\|\int_0^{\tilde{y}}\partial^{m-j+1}_x\tilde{b}ds\|_{L^2_xL^\infty_y}\\
\leq&C\langle t\rangle^{1/4}\|\theta_{\alpha}\partial^{m-j+1}_x\tilde{u}\|_{L^2}+C\langle t\rangle^{1/4}\|\theta_{\alpha}\partial^{m-j+1}_x\tilde{b}\|_{L^2},
\end{align*}
where we have used $\int_0^\infty|\partial_yu_s(t,y)|dy<C$ by the assumption (H). Note that
\begin{align*}
\|\theta_{\alpha}\partial^{j}_x\partial_y\tilde{u}\|_{L^\infty_xL^2_y}\leq C\|\theta_{\alpha}\partial^{j}_x\partial_y\tilde{u}\|^{1/2}_{L^2}\|\theta_{\alpha}\partial^{j+1}_x\partial_y\tilde{u}\|^{1/2}_{L^2}.
\end{align*}
For $[m/2]+1\leq j\leq m$,
\begin{align*}
\|\partial^{m-j}_xv\|_{L^\infty_{xy}}\leq & \|\int_0^y\partial^{m-j+1}_x\tilde{u}d\tilde{y}\|_{L^\infty_{xy}}+\|\int_0^y\partial_yu_s(\int_0^{\tilde{y}}\partial^{m-j+1}_x\tilde{b}ds)d\tilde{y}\|_{L^\infty_{xy}}\\
\leq & \|\int_0^y\partial^{m-j+1}_x\tilde{u}d\tilde{y}\|_{L^\infty_{xy}}+\|\int_0^y\partial_yu_s d\tilde{y}\|_{L^\infty_{y}}\|\int_0^{\tilde{y}}\partial^{m-j+1}_x\tilde{b}ds\|_{L^\infty_{xy}}\\
\leq & C\langle t\rangle^{1/4}\|\theta_{\alpha}\partial^{m-j+1}_x\tilde{u}\|_{L^2_yL^\infty_x}+C\langle t\rangle^{1/4}\|\theta_{\alpha}\partial^{m-j+1}_x\tilde{b}\|_{L^2_yL^\infty_x}\\
\leq &C\langle t\rangle^{1/4}\|\theta_{\alpha}\partial^{m-j+1}_x\tilde{u}\|^{1/2}_{L^2}\|\theta_{\alpha}\partial^{m-j+2}_x\tilde{u}\|^{1/2}_{L^2}+C\langle t\rangle^{1/4}\|\theta_{\alpha}\partial^{m-j+1}_x\tilde{b}\|^{1/2}_{L^2}\|\theta_{\alpha}\partial^{m-j+2}_x\tilde{b}\|^{1/2}_{L^2}.
\end{align*}
Hence,
\begin{align}
\label{NE2}
\frac{|R_2|\tau^mM_m}{\|\theta_{\alpha}\partial_x^m\tilde{u}\|_{L^2}}\leq& \frac{C}{(\tau(t))^{1/2}}\left\{\sum_{j=0}^{[m/2]}(\langle t\rangle^{1/4}Y_{m-j+1}+\langle t\rangle^{1/4}\bar{Y}_{m-j+1})D^{1/2}_{j}D^{1/2}_{j+1} \right.\\
&\left.+\sum_{j=[m/2]+1}^{m}(\langle t\rangle^{1/4}Y_{m-j+1}^{1/2}Y_{m-j+2}^{1/2}+\langle t\rangle^{1/4}\bar{Y}_{m-j+1}^{1/2}\bar{Y}_{m-j+2}^{1/2})D_{j}\right\}.\nonumber
\end{align}
Recall  $b=\tilde{b}$ so that
\begin{align*}
R_3\triangleq\sum_{j=0}^m(\begin{array}{ll}m\\j\end{array})\int_{\mathbb{H}}\partial^{m-j}_x\tilde{b}\partial_x^{j+1}\tilde{b}\theta_{\alpha}^2\partial_x^m\tilde{u}dxdy,
\end{align*}
and
\begin{align*}
|R_3|\leq& \sum_{j=0}^{[m/2]}(\begin{array}{ll}m\\j\end{array})\|\partial^{m-j}_x\tilde{b}\|_{L^2_xL^\infty_y}\|\theta_{\alpha}\partial^{j+1}_x\tilde{b}\|_{L^\infty_xL^2_y}\|\theta_{\alpha}\partial^{m}_x\tilde{u}\|_{L^2}\\
+&\sum_{j=[m/2]+1}^{m}(\begin{array}{ll}m\\j\end{array})\|\partial^{m-j}_x\tilde{b}\|_{L^\infty_{xy}}\|\theta_{\alpha}\partial^{j+1}_x\tilde{b}\|_{L^2}\|\theta_{\alpha}\partial^{m}_x\tilde{u}\|_{L^2}.
\end{align*}
For $0\leq j\leq [m/2]$,
\begin{align*}
\|\partial^{m-j}_x\tilde{b}\|_{L^2_xL^\infty_y}
\leq C\|\theta_{\alpha}\partial^{m-j}_x\tilde{b}\|^{1/2}_{L^2}\|\partial^{m-j}_x\partial_y\tilde{b}\|^{1/2}_{L^2},
\end{align*}
and
\begin{align*}
\|\partial^{j+1}_x\tilde{b}\|_{L^\infty_xL^2_y}\leq C\|\partial^{j+1}_x\tilde{b}\|^{1/2}_{L^2}\|\partial^{j+2}_x\tilde{b}\|^{1/2}_{L^2}.
\end{align*}
For $[m/2]+1\leq j\leq m$,
\begin{align*}
\|\partial^{m-j}_x\tilde{b}\|_{L^\infty_{xy}}
\leq C\|\partial^{m-j}_x\tilde{b}\|^{1/4}_{L^2}\|\partial^{m-j}_x\partial_y\tilde{b}\|^{1/4}_{L^2}\|\partial^{m-j+1}_x\tilde{b}\|^{1/4}_{L^2}\|\partial^{m-j+1}_x\partial_y\tilde{b}\|^{1/4}_{L^2}.
\end{align*}
Therefore,
\begin{align}
\label{NE3}
\frac{|R_3|\tau^mM_m}{\|\theta_{\alpha}\partial_x^m\tilde{u}\|_{L^2}}\leq& \frac{C}{(\tau(t))^{1/2}}\left\{\sum_{j=0}^{[m/2]}\bar{X}_{m-j}^{1/2}\bar{D}_{m-j}^{1/2}\bar{Y}_{j+1}^{1/2}\bar{Y}_{j+2}^{1/2}\right.\\
&\left.+\sum_{j=[m/2]+1}^{m}\bar{X}_{m-j}^{1/4}\bar{X}_{m-j+1}^{1/4}\bar{D}_{m-j}^{1/4}\bar{D}_{m-j+1}^{1/4}\bar{Y}_{j+1}\right\}.\nonumber
\end{align}
 Note that
\begin{align*}
\int_{\mathbb{H}}\partial_x^m(g\partial_y\tilde{b})\theta_{\alpha}^2\partial_x^m\tilde{u}dxdy
=\sum_{j=0}^m(\begin{array}{ll}m\\j\end{array})\int_{\mathbb{H}}\partial^{m-j}_xg\partial_x^{j}\partial_y\tilde{b}\theta_{\alpha}^2\partial_x^m\tilde{u}dxdy
\triangleq R_4
\end{align*}
and
\begin{align*}
|R_4|\leq& \sum_{j=0}^{[m/2]}(\begin{array}{ll}m\\j\end{array})\|\partial^{m-j}_xg\|_{L^2_xL^\infty_y}\|\theta_{\alpha}\partial^{j}_x\partial_y\tilde{b}\|_{L^\infty_xL^2_y}\|\theta_{\alpha}\partial^{m}_x\tilde{u}\|_{L^2}\\
+&\sum_{j=[m/2]+1}^{m}(\begin{array}{ll}m\\j\end{array})\|\partial^{m-j}_xg\|_{L^\infty_{xy}}\|\theta_{\alpha}\partial^{j}_x\partial_y\tilde{b}\|_{L^2}\|\theta_{\alpha}\partial^{m}_x\tilde{u}\|_{L^2}.
\end{align*}
For $0\leq j\leq [m/2]$,
\begin{align*}
\|\partial^{m-j}_xg\|_{L^2_xL^\infty_y}\leq  C\langle t\rangle^{1/4}\|\theta_{\alpha}\partial^{m-j+1}_x\tilde{b}\|_{L^2},
\end{align*}
and
\begin{align*}
\|\theta_{\alpha}\partial^{j}_x\partial_y\tilde{b}\|_{L^\infty_xL^2_y}\leq C\|\theta_{\alpha}\partial^{j}_x\partial_y\tilde{b}\|^{1/2}_{L^2}\|\theta_{\alpha}\partial^{j+1}_x\partial_y\tilde{b}\|^{1/2}_{L^2}.
\end{align*}
For $[m/2]+1\leq j\leq m$,
\begin{align*}
\|\partial^{m-j}_xg\|_{L^\infty_{xy}}\leq &
C\langle t\rangle^{1/4}\|\theta_{\alpha}\partial^{m-j+1}_x\tilde{b}\|_{L^2_yL^\infty_x}\\
\leq & C\langle t\rangle^{1/4}\|\theta_{\alpha}\partial^{m-j+1}_x\tilde{b}\|^{1/2}_{L^2}\|\theta_{\alpha}\partial^{m-j+2}_x\tilde{b}\|^{1/2}_{L^2}.
\end{align*}
As a consequence, we have
\begin{align}
\label{NE4}
\frac{|R_4|\tau^mM_m}{\|\theta_{\alpha}\partial_x^m\tilde{u}\|_{L^2}}\leq  \frac{C}{(\tau(t))^{1/2}}\left\{\sum_{j=0}^{[m/2]}\langle t\rangle^{1/4}\bar{Y}_{m-j+1}\bar{D}^{1/2}_{j}\bar{D}^{1/2}_{j+1}
+\sum_{j=[m/2]+1}^{m}\langle t\rangle^{1/4}\bar{Y}^{1/2}_{m-j+1}\bar{Y}^{1/2}_{m-j+2}\bar{D}_{j}\right\}.
\end{align}
And
\begin{align*}
|\int_\mathbb{H}\partial_y^2u_s\partial_x^m\tilde{b}\theta_{\alpha}^2\partial_x^m\tilde{u}dxdy|
\leq&\|\partial_y^2u_s\|_{L^\infty_y}\|\theta_{\alpha}\partial^{m}_x\tilde{b}\|_{L^2}\|\theta_{\alpha}\partial^{m}_x\tilde{u}\|_{L^2}\\
\leq&
C\langle t\rangle^{-1}\|\theta_{\alpha}\partial^{m}_x\tilde{b}\|_{L^2}\|\theta_{\alpha}\partial^{m}_x\tilde{u}\|_{L^2},
\end{align*}
that is,
\begin{align}
\label{NE5}
\frac{|\int_\mathbb{H}\partial_y^2u_s\partial_x^m\tilde{b}\theta_{\alpha}^2\partial_x^m\tilde{u}dxdy|}{\|\theta_{\alpha}\partial_x^m\tilde{u}\|_{L^2}}\leq C\langle t\rangle^{-1}\|\theta_{\alpha}\partial^{m}_x\tilde{b}\|_{L^2},
\end{align}
where we  have used $\|\partial_y^2u_s\|_{L_y^\infty}\leq \frac{C}{\langle t\rangle}$ by the assumption (H). We now consider
\begin{align*}
R_5\triangleq\sum_{j=0}^m(\begin{array}{ll}m\\j\end{array})\int_{H}\partial^{m-j}_xv \partial_y^2u_s   \partial_x^{j}\psi\theta_{\alpha}^2\partial_x^m\tilde{u}dxdy.
\end{align*}
Note that
\begin{align*}
|R_5|\leq& \sum_{j=0}^{[m/2]}(\begin{array}{ll}m\\j\end{array})\|\partial^{m-j}_xv\|_{L^2_xL^\infty_y}\|\theta_{\alpha}\partial_y^2u_s\|_{L^2_y}  \|\partial^{j}_x\psi\|_{L^\infty_{xy}}\|\theta_{\alpha}\partial^{m}_x\tilde{u}\|_{L^2}\\
+&\sum_{j=[m/2]+1}^{m}(\begin{array}{ll}m\\j\end{array})\|\partial^{m-j}_xv\|_{L^\infty_{xy}}\|\theta_{\alpha}\partial_y^2u_s\|_{L^2_y}  \|\partial^{j}_x\psi\|_{L^2_xL^\infty_y}\|\theta_{\alpha}\partial^{m}_x\tilde{u}\|_{L^2}.
\end{align*}
For $0\leq j\leq [m/2]$, we have
\begin{align*}
\|\partial^{m-j}_xv\|_{L^2_xL^\infty_y}
\leq C\langle t\rangle^{1/4}\|\theta_{\alpha}\partial^{m-j+1}_x\tilde{u}\|_{L^2} + C\langle t\rangle^{1/4}\|\theta_{\alpha}\partial^{m-j+1}_x\tilde{b}\|_{L^2},
\end{align*}
and
\begin{align*}
\|\theta_{\alpha}\partial_y^2u_s\|_{L^2_y}\leq \frac{C}{\langle t\rangle^{3/4}},
\end{align*}
provided that $\alpha<1$ by the  assumption (H). And
\begin{align*}
\|\partial^{j}_x\psi\|_{L^\infty_{xy}}\leq C\langle t\rangle^{1/4} \|\theta_{\alpha}\partial^{j}_x\tilde{b}\|^{1/2}_{L^2}\|\theta_{\alpha}\partial^{j+1}_x\tilde{b}\|^{1/2}_{L^2}.
\end{align*}
For $[m/2]+1\leq j\leq m$, we have
\begin{align*}
\|\partial^{m-j}_xv\|_{L^\infty_{xy}}
\leq& C\langle t\rangle^{1/4} \|\theta_{\alpha}\partial^{m-j+1}_x\tilde{u}\|^{1/2}_{L^2}\|\theta_{\alpha}\partial^{m-j+2}_x\tilde{u}\|^{1/2}_{L^2}\\
&+C\langle t\rangle^{1/4} \|\theta_{\alpha}\partial^{m-j+1}_x\tilde{b}\|^{1/2}_{L^2}\|\theta_{\alpha}\partial^{m-j+2}_x\tilde{b}\|^{1/2}_{L^2}.
\end{align*}
And
\begin{align*}
\|\partial^{j}_x\psi\|_{L^2_xL^\infty_y}\leq C\langle t\rangle^{1/4}\|\theta_{\alpha}\partial_x^j\tilde{b}\|_{L^2}.
\end{align*}
Hence,
\begin{align}
\label{NE6}
\frac{|R_5|\tau^mM_m}{\|\theta_{\alpha}\partial_x^m\tilde{u}\|_{L^2}}\leq & \frac{C}{(\tau(t))^{1/2}}\left\{\sum_{j=0}^{[m/2]}(\langle t\rangle^{-1/4}Y_{m-j+1}+\langle t\rangle^{-1/4}\bar{Y}_{m-j+1})\bar{X}^{1/2}_{j}\bar{X}^{1/2}_{j+1} \right.\\
&\left.+\sum_{j=[m/2]+1}^{m}(\langle t\rangle^{-1/4}Y_{m-j+1}^{1/2}Y_{m-j+2}^{1/2}+\langle t\rangle^{-1/4}\bar{Y}_{m-j+1}^{1/2}\bar{Y}_{m-j+2}^{1/2})\bar{X}_{j}\right\}.\nonumber
\end{align}
Combining the estimates (\ref{3.6})-(\ref{NE6})  and summing over $m\geq 0$ give
\begin{align}
\label{NEU}
&\frac{d}{dt}\|\tilde{u}\|_{X_{\tau,\alpha}}+\sum_{m\geq 0}\tau^mM_m\frac{\|\theta_{\alpha}\partial_x^m\partial_y\tilde{u}\|_{L^2}^2}{\|\theta_\alpha\partial_x^m\tilde{u}\|_{L^2}}+\frac{\alpha(1-2\alpha)}{4\langle t\rangle}\sum_{m\geq 0}\tau^mM_m\frac{\|\theta_{\alpha}z\partial_x^m\tilde{u}\|^2_{L^2}}{\|\theta_\alpha\partial_x^m\tilde{u}\|_{L^2}}\nonumber\\
&-\frac{\alpha}{2\langle t\rangle}\|\tilde{u}\|_{X_{\tau,\alpha}}-\frac{C}{\langle t\rangle}\|\tilde{b}\|_{X_{\tau,\alpha}}
\leq\dot{\tau}(t)\|\tilde{u}\|_{Y_{\tau,\alpha}}\\&+\frac{C_0}{(\tau(t))^{1/2}}\left(\langle t\rangle^{-1/4}(\|\tilde{u}\|_{X_{\tau,\alpha}}+\|\tilde{b}\|_{X_{\tau,\alpha}})+\langle t\rangle^{1/4}(\|\tilde{u}\|_{D_{\tau,\alpha}}+\|\tilde{b}\|_{D_{\tau,\alpha}})\right)(\|\tilde{u}\|_{Y_{\tau,\alpha}}+\|\tilde{b}\|_{Y_{\tau,\alpha}}),\nonumber
\end{align}
where we have used the fact that for any positive sequences $\{a_j\}_{j\geq 0}$ and $\{b_j\}_{j\geq 0}$,
\begin{align*}
\sum_{m\geq 0}\sum_{j=0}^{m}a_jb_{m-j}\leq \sum_{j\geq 0}a_j\sum_{j\geq 0}b_j.
\end{align*}
Choosing $\alpha\leq1/2$ in (\ref{NEU}) yields
\begin{align}
\label{NEUU}
&\frac{d}{dt}\|\tilde{u}\|_{X_{\tau,\alpha}}+\sum_{m\geq 0}\tau^mM_m\frac{\|\theta_{\alpha}\partial_x^m\partial_y\tilde{u}\|_{L^2}^2}{\|\theta_\alpha\partial_x^m\tilde{u}\|_{L^2}}-\frac{\alpha}{2\langle t\rangle}\|\tilde{u}\|_{X_{\tau,\alpha}}-\frac{C}{\langle t\rangle}\|\tilde{b}\|_{X_{\tau,\alpha}}\nonumber\\
\leq& \dot{\tau}(t)\|\tilde{u}\|_{Y_{\tau,\alpha}}+\frac{C_0}{(\tau(t))^{1/2}}\left(\langle t\rangle^{-1/4}(\|\tilde{u}\|_{X_{\tau,\alpha}}+\|\tilde{b}\|_{X_{\tau,\alpha}})+\langle t\rangle^{1/4}(\|\tilde{u}\|_{D_{\tau,\alpha}}+\|\tilde{b}\|_{D_{\tau,\alpha}})\right)\nonumber\\
&\times(\|\tilde{u}\|_{Y_{\tau,\alpha}}+\|\tilde{b}\|_{Y_{\tau,\alpha}}).
\end{align}
\subsection{A priori estimate on  magnetic field}
Similarly,
for $m\geq 0$,  by applying the tangential derivative operator $\partial_x^m$ to $(\ref{3.4})_2$ and  multiplying it by $\theta_{\alpha}^2\partial_x^m\tilde{b}$, the integration  over  $\mathbb{H}$ gives
\begin{align}
\label{3.8}
\int_{\mathbb{H}}\partial_x^m(\partial_t\tilde{b}-\partial_y^2\tilde{b}-(1+b)\partial_x\tilde{u}-g\partial_y\tilde{u}+
(u_s+u)\partial_x\tilde{b}+v\partial_y\tilde{b}-g\partial_y^2u_s\psi)\theta_{\alpha}^2\partial_x^m\tilde{b}dxdy=0.
\end{align}
We now estimate $(\ref{3.8})$ term by term as follows. Firstly,
\begin{align}
\label{3.19}
&\int_\mathbb{H}\partial_t\partial_x^m\tilde{b}\theta_{\alpha}^2\partial_x^m\tilde{b}dxdy
= \frac12\frac{d}{dt}\int_\mathbb{H}(\partial_x^m\tilde{b})^2\theta_{\alpha}^2dxdy-\int_\mathbb{H}(\partial_x^m\tilde{b})^2\theta_{\alpha}\frac{d}{dt}\theta_{\alpha}dxdy\\
= &\frac12\frac{d}{dt}\|\theta_{\alpha}\partial_x^m\tilde{b}\|^2_{L^2}+\frac{\alpha}{4\langle t\rangle}\|\theta_{\alpha}z\partial_x^m\tilde{b}\|^2_{L^2}.\nonumber
\end{align}
And
\begin{align*}
-\int_\mathbb{H}\partial_y^2\partial_x^m\tilde{b}\theta_{\alpha}^2\partial_x^m\tilde{b}dxdy
=\|\theta_{\alpha}\partial_x^m\partial_y\tilde{b}\|_{L^2}^2+\int_\mathbb{H}\partial_y\partial_x^m\tilde{b}\partial_y(\theta_{\alpha}^2)\partial_x^m\tilde{b}dxdy,
\end{align*}
where we have used the boundary condition $\partial_y\partial_x^m\tilde{b}|_{y=0}=0$. Moreover,
\begin{align*}
&\int_\mathbb{H}\partial_y\partial_x^m\tilde{b}\partial_y(\theta_{\alpha}^2)\partial_x^m\tilde{b}dxdy
=-\frac12\int_\mathbb{H}(\partial_x^m\tilde{b})^2\partial_y^2(\theta_{\alpha}^2)dxdy\\
=&-\frac{\alpha}{2}\frac{1}{\langle t\rangle}\|\theta_{\alpha}\partial_x^m\tilde{b}\|^2_{L^2}-\frac{\alpha^2}{2}\frac{1}{\langle t\rangle}\|\theta_{\alpha}z\partial_x^m\tilde{b}\|^2_{L^2}.\nonumber
\end{align*}
Hence,
\begin{align}
\label{3.7}
-\int_\mathbb{H}\partial_y^2\partial_x^m\tilde{b}\theta_{\alpha}^2\partial_x^m\tilde{b}dxdy
= \|\theta_{\alpha}\partial_x^m\partial_y\tilde{b}\|_{L^2}^2-\frac{\alpha}{2}\frac{1}{\langle t\rangle}\|\theta_{\alpha}\partial_x^m\tilde{b}\|^2_{L^2}-\frac{\alpha^2}{2}\frac{1}{\langle t\rangle}\|\theta_{\alpha}z\partial_x^m\tilde{b}\|^2_{L^2}.
\end{align}
Similar to Subsection 2.1, the nonlinear terms can be estimated as follows.
Firstly,
\begin{align*}
\int_{\mathbb{H}}\partial_x^m((1+b)\partial_x\tilde{u})\theta_{\alpha}^2\partial_x^m\tilde{b}dxdy
=\sum_{j=0}^m(\begin{array}{ll}m\\j\end{array})\int_{H}\partial^{m-j}_x\tilde{b}\partial_x^{j+1}\tilde{u}\theta_{\alpha}^2\partial_x^m\tilde{b}dxdy
\triangleq R_6,
\end{align*}
and
\begin{align*}
|R_6|\leq& \sum_{j=0}^{[m/2]}(\begin{array}{ll}m\\j\end{array})\|\partial^{m-j}_x\tilde{b}\|_{L^2_xL^\infty_y}\|\theta_{\alpha}\partial^{j+1}_x\tilde{u}\|_{L^\infty_xL^2_y}\|\theta_{\alpha}\partial^{m}_x\tilde{b}\|_{L^2}\\
+&\sum_{j=[m/2]+1}^{m}(\begin{array}{ll}m\\j\end{array})\|\partial^{m-j}_x\tilde{b}\|_{L^\infty_{xy}}\|\theta_{\alpha}\partial^{j+1}_x\tilde{u}\|_{L^2}\|\theta_{\alpha}\partial^{m}_x\tilde{b}\|_{L^2}.
\end{align*}
For $0\leq j\leq [m/2]$,
\begin{align*}
\|\partial^{m-j}_x\tilde{b}\|_{L^2_xL^\infty_y}
\leq C\|\theta_{\alpha}\partial^{m-j}_x\tilde{b}\|^{1/2}_{L^2}\|\theta_{\alpha}\partial^{m-j}_x\partial_y\tilde{b}\|^{1/2}_{L^2},
\end{align*}
and
\begin{align*}
\|\theta_{\alpha}\partial^{j+1}_x\tilde{u}\|_{L^\infty_xL^2_y}\leq C\|\theta_{\alpha}\partial^{j+1}_x\tilde{u}\|^{1/2}_{L^2}\|\theta_{\alpha}\partial^{j+2}_x\tilde{u}\|^{1/2}_{L^2}.
\end{align*}
For $[m/2]+1\leq j\leq m$,
\begin{align*}
\|\partial^{m-j}_x\tilde{b}\|_{L^\infty_{xy}}\leq
C\|\theta_{\alpha}\partial^{m-j}_x\tilde{b}\|^{1/4}_{L^2}\|\theta_{\alpha}\partial^{m-j}_x\partial_y\tilde{b}\|^{1/4}_{L^2}
\|\theta_{\alpha}\partial^{m-j+1}_x\tilde{b}\|^{1/4}_{L^2}\|\theta_{\alpha}\partial^{m-j+1}_x\partial_y\tilde{b}\|^{1/4}_{L^2}.
\end{align*}
Hence,
\begin{align}
\label{NE7}
\frac{|R_6|\tau^mM_m}{\|\theta_{\alpha}\partial_x^m\tilde{b}\|_{L^2}}\leq \frac{C}{(\tau(t))^{1/2}}\left\{\sum_{j=0}^{[m/2]}\bar{X}_{m-j}^{1/2}\bar{D}_{m-j}^{1/2}Y_{j+1}^{1/2}Y_{j+2}^{1/2}+\sum_{j=[m/2]+1}^{m}\bar{X}^{1/4}_{m-j}\bar{X}^{1/4}_{m-j+1}\bar{D}^{1/4}_{m-j}\bar{D}^{1/4}_{m-j+1}Y_{j+1}\right\}.
\end{align}
Moreover,
\begin{align*}
\int_{\mathbb{H}}\partial_x^m(g\partial_y\tilde{u})\theta_{\alpha}^2\partial_x^m\tilde{b}dxdy
=\sum_{j=0}^m(\begin{array}{ll}m\\j\end{array})\int_{\mathbb{H}}\partial^{m-j}_xg\partial_x^{j}\partial_y\tilde{u}\theta_{\alpha}^2\partial_x^m\tilde{b}dxdy
\triangleq R_7,
\end{align*}
and
\begin{align*}
|R_7|\leq& \sum_{j=0}^{[m/2]}(\begin{array}{ll}m\\j\end{array})\|\partial^{m-j}_xg\|_{L^2_xL^\infty_y}\|\theta_{\alpha}\partial^{j}_x\partial_y\tilde{u}\|_{L^\infty_xL^2_y}\|\theta_{\alpha}\partial^{m}_x\tilde{b}\|_{L^2}\\
+&\sum_{j=[m/2]+1}^{m}(\begin{array}{ll}m\\j\end{array})\|\partial^{m-j}_xg\|_{L^\infty_{xy}}\|\theta_{\alpha}\partial^{j}_x\partial_y\tilde{u}\|_{L^2}\|\theta_{\alpha}\partial^{m}_x\tilde{b}\|_{L^2}.
\end{align*}
For $0\leq j\leq [m/2]$,
\begin{align*}
\|\partial^{m-j}_xg\|_{L^2_xL^\infty_y}\leq C\langle t\rangle^{1/4}\|\theta_{\alpha}\partial^{m-j+1}_x\tilde{b}\|_{L^2},
\end{align*}
and
\begin{align*}
\|\theta_{\alpha}\partial^{j}_x\partial_y\tilde{u}\|_{L^\infty_xL^2_y}\leq C\|\theta_{\alpha}\partial^{j}_x\partial_y\tilde{u}\|^{1/2}_{L^2}\|\theta_{\alpha}\partial^{j+1}_x\partial_y\tilde{u}\|^{1/2}_{L^2}.
\end{align*}
For $[m/2]+1\leq j\leq m$,
\begin{align*}
\|\partial^{m-j}_xg\|_{L^\infty_{xy}}\leq C \langle t\rangle^{1/4}\|\theta_{\alpha}\partial^{m-j+1}_x\tilde{b}\|^{1/2}_{L^2}\|\theta_{\alpha}\partial^{m-j+2}_x\tilde{b}\|^{1/2}_{L^2}.
\end{align*}
Therefore,
\begin{align}
\label{NE8}
\frac{|R_7|\tau^mM_m}{\|\theta_{\alpha}\partial_x^m\tilde{b}\|_{L^2}}\leq \frac{C}{(\tau(t))^{1/2}}\left\{\sum_{j=0}^{[m/2]}\langle t\rangle^{1/4}\bar{Y}_{m-j+1}D_{j}^{1/2}D_{j+1}^{1/2}+\sum_{j=[m/2]+1}^{m}\langle t\rangle^{1/4}\bar{Y}_{m-j+1}^{1/2}\bar{Y}_{m-j+2}^{1/2}D_{j}\right\}.
\end{align}
Denote
\begin{align*}
R_8\triangleq\sum_{j=0}^m(\begin{array}{ll}m\\j\end{array})\int_{\mathbb{H}}\partial^{m-j}_xu\partial_x^{j+1}\tilde{b}\theta_{\alpha}^2\partial_x^m\tilde{b}dxdy,
\end{align*}
then
\begin{align*}
|R_8|\leq& \sum_{j=0}^{[m/2]}(\begin{array}{ll}m\\j\end{array})\|\partial^{m-j}_xu\|_{L^2_xL^\infty_y}\|\theta_{\alpha}\partial^{j+1}_x\tilde{b}\|_{L^\infty_xL^2_y}\|\theta_{\alpha}\partial^{m}_x\tilde{b}\|_{L^2}\\
+&\sum_{j=[m/2]+1}^{m}(\begin{array}{ll}m\\j\end{array})\|\partial^{m-j}_xu\|_{L^\infty_{xy}}\|\theta_{\alpha}\partial^{j+1}_x\tilde{b}\|_{L^2}\|\theta_{\alpha}\partial^{m}_x\tilde{b}\|_{L^2}.
\end{align*}
Similar to the estimation on $R_1$, we can obtain
\begin{align}
\label{NE9}
\frac{|R_8|\tau^mM_m}{\|\theta_{\alpha}\partial_x^m\tilde{b}\|_{L^2}}
\leq& \frac{C}{(\tau(t))^{1/2}}\left\{\sum_{j=0}^{[m/2]}(X_{m-j}^{1/2}D_{m-j}^{1/2}+\langle t\rangle^{-1/4}\bar{X}_{m-j})\bar{Y}_{j+1}^{1/2}\bar{Y}_{j+2}^{1/2}\right.\\
&\left.+\sum_{j=[m/2]+1}^{m}(X_{m-j}^{1/4}X_{m-j+1}^{1/4}D_{m-j}^{1/4}D_{m-j+1}^{1/4}+\langle t\rangle^{-1/4}\bar{X}^{1/2}_{m-j}\bar{X}^{1/2}_{m-j+1})\bar{Y}_{j+1}\right\}.\nonumber
\end{align}
And
\begin{align*}
\int_{\mathbb{H}}\partial_x^m(v\partial_y\tilde{b})\theta_{\alpha}^2\partial_x^m\tilde{b}dxdy
=\sum_{j=0}^m(\begin{array}{ll}m\\j\end{array})\int_{\mathbb{H}}\partial^{m-j}_xv\partial_x^{j}\partial_y\tilde{b}\theta_{\alpha}^2\partial_x^m\tilde{b}dxdy
\triangleq R_9.
\end{align*}
Thus
\begin{align*}
|R_9|\leq& \sum_{j=0}^{[m/2]}(\begin{array}{ll}m\\j\end{array})\|\partial^{m-j}_xv\|_{L^2_xL^\infty_y}\|\theta_{\alpha}\partial^{j}_x\partial_y\tilde{b}\|_{L^\infty_xL^2_y}\|\theta_{\alpha}\partial^{m}_x\tilde{b}\|_{L^2}\\
+&\sum_{j=[m/2]+1}^{m}(\begin{array}{ll}m\\j\end{array})\|\partial^{m-j}_xv\|_{L^\infty_{xy}}\|\theta_{\alpha}\partial^{j}_x\partial_y\tilde{b}\|_{L^2}\|\theta_{\alpha}\partial^{m}_x\tilde{b}\|_{L^2}.
\end{align*}
Similar to  the estimation on  $R_2$, we have
\begin{align}
\label{NE10}
\frac{|R_9|\tau^mM_m}{\|\theta_{\alpha}\partial_x^m\tilde{b}\|_{L^2}}
\leq& \frac{C}{(\tau(t))^{1/2}}\left\{\sum_{j=0}^{[m/2]}(\langle t\rangle^{1/4}Y_{m-j+1}+\langle t\rangle^{1/4}\bar{Y}_{m-j+1})\bar{D}_{j}^{1/2}\bar{D}_{j+1}^{1/2}\right.\\
&\left.+\sum_{j=[m/2]+1}^{m}(\langle t\rangle^{1/4}Y_{m-j+1}^{1/2}Y_{m-j+2}^{1/2}+\langle t\rangle^{1/4}\bar{Y}_{m-j+1}^{1/2}\bar{Y}_{m-j+2}^{1/2})\bar{D}_{j}\right\}.\nonumber
\end{align}
Denote
\begin{align*}
R_{10} \triangleq  \int_\mathbb{H}\partial_x^m(g\partial_y^2u_s\psi)\theta_{\alpha}^2\partial_x^m\tilde{b}dxdy
=\sum_{j=0}^m\int_\mathbb{H}\partial_x^{m-j}g\partial_y^2u_s\partial_x^j\psi\theta_{\alpha}^2\partial_x^m\tilde{b}dxdy.
\end{align*}
Then
\begin{align*}
|R_{10}|\leq& \sum_{j=0}^{[m/2]}(\begin{array}{ll}m\\j\end{array})\|\partial^{m-j}_xg\|_{L^2_xL^\infty_y}\|\theta_{\alpha}\partial_y^2u_s\|_{L^2_y}  \|\partial^{j}_x\psi\|_{L^\infty_{xy}}\|\theta_{\alpha}\partial^{m}_x\tilde{b}\|_{L^2}\\
+&\sum_{j=[m/2]+1}^{m}(\begin{array}{ll}m\\j\end{array})\|\partial^{m-j}_xg\|_{L^\infty_{xy}}\|\theta_{\alpha}\partial_y^2u_s\|_{L^2_y}  \|\partial^{j}_x\psi\|_{L^2_xL^\infty_y}\|\theta_{\alpha}\partial^{m}_x\tilde{b}\|_{L^2}.
\end{align*}
For $0\leq j\leq [m/2]$,
\begin{align*}
\|\partial^{m-j}_xg\|_{L^2_xL^\infty_y}
\leq C\langle t\rangle^{1/4}\|\theta_{\alpha}\partial^{m-j+1}_x\tilde{b}\|_{L^2},
\end{align*}
and
\begin{align*}
\|\theta_{\alpha}\partial_y^2u_s\|_{L^2_y}\leq \frac{C}{\langle t\rangle^{3/4}},
\end{align*}
provided that $\alpha<1$ by the assumption (H). Moreover,
\begin{align*}
\|\partial^{j}_x\psi\|_{L^\infty_{xy}}\leq C\langle t\rangle^{1/4} \|\theta_{\alpha}\partial^{j}_x\tilde{b}\|^{1/2}_{L^2}\|\theta_{\alpha}\partial^{j+1}_x\tilde{b}\|^{1/2}_{L^2}.
\end{align*}
For $[m/2]+1\leq j\leq m$,
\begin{align*}
\|\partial^{m-j}_xg\|_{L^\infty_{xy}}
\leq C\langle t\rangle^{1/4} \|\theta_{\alpha}\partial^{m-j+1}_x\tilde{b}\|^{1/2}_{L^2}\|\theta_{\alpha}\partial^{m-j+2}_x\tilde{b}\|^{1/2}_{L^2},
\end{align*}
and
\begin{align*}
\|\partial^{j}_x\psi\|_{L^2_xL^\infty_y}\leq C\langle t\rangle^{1/4}\|\theta_\alpha\partial_x^j\tilde{b}\|_{L^2}.
\end{align*}
Consequently,
\begin{align}
\label{NE11}
\frac{|R_{10}|\tau^mM_m}{\|\theta_{\alpha}\partial_x^m\tilde{b}\|_{L^2}}
\leq& \frac{C}{(\tau(t))^{1/2}}\left\{\sum_{j=0}^{[m/2]}\langle t\rangle^{-1/4}\bar{Y}_{m-j+1}\bar{X}_j^{1/2}\bar{X}_{j+1}^{1/2}\right.\\
&\left.+\sum_{j=[m/2]+1}^{m}\langle t\rangle^{-1/4}\bar{Y}^{1/2}_{m-j+1}\bar{Y}^{1/2}_{m-j+2}\bar{X}_j\right\}.\nonumber
\end{align}
From the estimates (\ref{3.19})-(\ref{NE11}), summing over $m\geq 0$ yields
\begin{align}
\label{NEB}
&\frac{d}{dt}\|\tilde{b}\|_{X_{\tau,\alpha}}+\sum_{m\geq 0}\tau^mM_m\frac{\|\theta_{\alpha}\partial_x^m\partial_y\tilde{b}\|_{L^2}^2}{\|\theta_\alpha\partial_x^m\tilde{b}\|_{L^2}}+\frac{\alpha(1-2\alpha)}{4\langle t\rangle}\sum_{m\geq 0}\tau^mM_m\frac{\|\theta_{\alpha}z\partial_x^m\tilde{b}\|^2_{L^2}}{\|\theta_\alpha\partial_x^m\tilde{b}\|_{L^2}}-\frac{\alpha}{2\langle t\rangle}\|\tilde{b}\|_{X_{\tau,\alpha}}\nonumber\\
\leq&\dot{\tau}(t)\|\tilde{b}\|_{Y_{\tau,\alpha}}+\frac{C_0}{(\tau(t))^{1/2}}\left(\langle t\rangle^{-1/4}(\|\tilde{u}\|_{X_{\tau,\alpha}}+\|\tilde{b}\|_{X_{\tau,\alpha}})+\langle t\rangle^{1/4}(\|\tilde{u}\|_{D_{\tau,\alpha}}+\|\tilde{b}\|_{D_{\tau,\alpha}})\right)\nonumber\\
&\times(\|\tilde{u}\|_{Y_{\tau,\alpha}}+\|\tilde{b}\|_{Y_{\tau,\alpha}}).
\end{align}
Similarly, by choosing $\alpha\leq1/2$, we have
\begin{align}
\label{NEBB}
&\frac{d}{dt}\|\tilde{b}\|_{X_{\tau,\alpha}}+\sum_{m\geq 0}\tau^mM_m\frac{\|\theta_{\alpha}\partial_x^m\partial_y\tilde{b}\|_{L^2}^2}{\|\theta_\alpha\partial_x^m\tilde{b}\|_{L^2}}-\frac{\alpha}{2\langle t\rangle}\|\tilde{b}\|_{X_{\tau,\alpha}}\nonumber\\
\leq& \dot{\tau}(t)\|\tilde{b}\|_{Y_{\tau,\alpha}}+\frac{C_0}{(\tau(t))^{1/2}}\left(\langle t\rangle^{-1/4}(\|\tilde{u}\|_{X_{\tau,\alpha}}+\|\tilde{b}\|_{X_{\tau,\alpha}})+\langle t\rangle^{1/4}(\|\tilde{u}\|_{D_{\tau,\alpha}}+\|\tilde{b}\|_{D_{\tau,\alpha}})\right)\nonumber\\
&\times(\|\tilde{u}\|_{Y_{\tau,\alpha}}+\|\tilde{b}\|_{Y_{\tau,\alpha}}).
\end{align}
\section{The Proof of Estimate of Lifespan in Theorem \ref{THM}}
By the uniform {\it a priori} estimates obtained in Section 2, we now estimate
the low bound on lifespan of the solution.  Consider
$(\ref{NEUU})+K\times(\ref{NEBB})$ with $K> 1$ to be determined later,
\begin{align}
\label{NEUB}
&\frac{d}{dt}(\|\tilde{u}\|_{X_{\tau,\alpha}}+K\|\tilde{b}\|_{X_{\tau,\alpha}})+\sum_{m\geq 0}\tau^mM_m(\frac{\|\theta_{\alpha}\partial_x^m\partial_y\tilde{u}\|_{L^2}^2}{\|\theta_\alpha\partial_x^m\tilde{u}\|_{L^2}}
+K\frac{\|\theta_{\alpha}\partial_x^m\partial_y\tilde{b}\|_{L^2}^2}{\|\theta_\alpha\partial_x^m\tilde{b}\|_{L^2}})\nonumber\\
&-\frac{\alpha}{2\langle t\rangle}\|\tilde{u}\|_{X_{\tau,\alpha}}-(C+\frac{K\alpha}{2})\frac{1}{\langle t\rangle}\|\tilde{b}\|_{X_{\tau,\alpha}}\\
\leq&\left(\dot{\tau}(t)+\frac{C_0(K+1)}{(\tau(t))^{1/2}}\left(\langle t\rangle^{-1/4}(\|\tilde{u}\|_{X_{\tau,\alpha}}+\|\tilde{b}\|_{X_{\tau,\alpha}})+\langle t\rangle^{1/4}(\|\tilde{u}\|_{D_{\tau,\alpha}}+\|\tilde{b}\|_{D_{\tau,\alpha}})\right)\right)\nonumber\\
&\times(\|\tilde{u}\|_{Y_{\tau,\alpha}}+K\|\tilde{b}\|_{Y_{\tau,\alpha}}).\nonumber
\end{align}
Choose the function $\tau(t)$ satisfies the following ODE.
\begin{align}
\label{ODE}
\frac{d}{dt}(\tau(t))^{3/2}+\frac{3C_0(K+1)}{2}\left(\langle t\rangle^{-1/4}(\|\tilde{u}\|_{X_{\tau,\alpha}}+\|\tilde{b}\|_{X_{\tau,\alpha}})+\langle t\rangle^{1/4}(\|\tilde{u}\|_{D_{\tau,\alpha}}+\|\tilde{b}\|_{D_{\tau,\alpha}})\right)=0.
\end{align}
From (\ref{NEUB}) and (\ref{ODE}), one has
\begin{align}
\label{NEUBU}
\frac{d}{dt}(\|\tilde{u}\|_{X_{\tau,\alpha}}+K\|\tilde{b}\|_{X_{\tau,\alpha}})&+\sum_{m\geq 0}\tau^mM_m\left(\frac{\|\theta_{\alpha}\partial_x^m\partial_y\tilde{u}\|_{L^2}^2}{\|\theta_\alpha\partial_x^m\tilde{u}\|_{L^2}}
+K\frac{\|\theta_{\alpha}\partial_x^m\partial_y\tilde{b}\|_{L^2}^2}{\|\theta_\alpha\partial_x^m\tilde{b}\|_{L^2}}\right)\nonumber\\
&-\frac{\alpha}{2\langle t\rangle}\|\tilde{u}\|_{X_{\tau,\alpha}}-(C+\frac{K\alpha}{2})\frac{1}{\langle t\rangle}\|\tilde{b}\|_{X_{\tau,\alpha}} \leq
0.
\end{align}
By  lemma \ref{LEM2.2}, we have
\begin{align}
\label{ESTU}
\sum_{m\geq 0}\tau^mM_m\frac{\|\theta_{\alpha}\partial_x^m\partial_y\tilde{u}\|_{L^2}^2}{\|\theta_\alpha\partial_x^m\tilde{u}\|_{L^2}}\geq \frac{\alpha^{1/2}\beta_1}{2\langle t\rangle^{1/2}}\|\tilde{u}\|_{D_{\tau,\alpha}}+\frac{\alpha(1-\beta_1)}{\langle t\rangle}\|\tilde{u}\|_{X_{\tau,\alpha}},
\end{align}
and
\begin{align}
\label{ESTB}
\sum_{m\geq 0}\tau^mM_m\frac{\|\theta_{\alpha}\partial_x^m\partial_y\tilde{b}\|_{L^2}^2}{\|\theta_\alpha\partial_x^m\tilde{b}\|_{L^2}}\geq \frac{\alpha^{1/2}\beta_2}{2\langle t\rangle^{1/2}}\|\tilde{b}\|_{D_{\tau,\alpha}}+\frac{\alpha(1-\beta_2)}{\langle t\rangle}\|\tilde{b}\|_{X_{\tau,\alpha}},
\end{align}
for $\beta_1, \beta_2\in(0,1/2)$.

From (\ref{NEUBU}), (\ref{ESTU}) and (\ref{ESTB}), it follows that
\begin{align}
\label{ESTT}
\frac{d}{dt}(\|\tilde{u}\|_{X_{\tau,\alpha}}+K\|\tilde{b}\|_{X_{\tau,\alpha}})&+\frac12(\alpha(1-2\beta_1))\frac{1}{\langle t\rangle}\|\tilde{u}\|_{X_{\tau,\alpha}}+
\frac12(\alpha(1-2\beta_2)-\frac{2C}{K})\frac{1}{\langle t\rangle}K\|\tilde{b}\|_{X_{\tau,\alpha}}\nonumber\\
&\qquad\qquad+\frac{\alpha^{1/2}\beta_1}{2\langle t\rangle^{1/2}}\|\tilde{u}\|_{D_{\tau,\alpha}}+\frac{K\alpha^{1/2}\beta_2}{2\langle t\rangle^{1/2}}\|\tilde{b}\|_{D_{\tau,\alpha}}\leq 0.
\end{align}
Choose
\begin{align*}
\alpha=\frac12-\delta,\quad \beta_1=\frac{\delta}{2},\quad \beta_2=\frac{\delta}{2},\quad K=\frac{4C}{\delta},
\end{align*}
where $0<\delta<1/4$ is sufficiently small to be determined later, then
\begin{align*}
\alpha(1-2\beta_1)=\frac12-\frac32\delta+\delta^2,
\end{align*}
and
\begin{align*}
\alpha(1-2\beta_2)-\frac{2C}{K}=\frac12-2\delta+\delta^2.
\end{align*}
Then, there exist small positive constants $\eta_1=\delta$ and $\eta_2=\frac{\delta}{8}$  such that
\begin{align}
\label{ESTTT}
\frac{d}{dt}(\|\tilde{u}\|_{X_{\tau,\alpha}}+K\|\tilde{b}\|_{X_{\tau,\alpha}})+\frac{1/4-\eta_1}{\langle t\rangle}\left(\|\tilde{u}\|_{X_{\tau,\alpha}}+
K\|\tilde{b}\|_{X_{\tau,\alpha}}\right)+\frac{\eta_2}{\langle t\rangle^{1/2}}(\|\tilde{u}\|_{D_{\tau,\alpha}}+K\|\tilde{b}\|_{D_{\tau,\alpha}})\leq 0.
\end{align}
It implies that
\begin{align}
\label{ESTTTT}
\frac{d}{dt}(\|\tilde{u}\|_{X_{\tau,\alpha}}+K\|\tilde{b}\|_{X_{\tau,\alpha}})\langle t\rangle^{1/4-\eta_1}&+\frac{1/4-\eta_1}{\langle t\rangle^{3/4+\eta_1}}\left(\|\tilde{u}\|_{X_{\tau,\alpha}}+
K\|\tilde{b}\|_{X_{\tau,\alpha}}\right)\nonumber\\
&+\frac{\eta_2}{\langle t\rangle^{1/4+\eta_1}}(\|\tilde{u}\|_{D_{\tau,\alpha}}+K\|\tilde{b}\|_{D_{\tau,\alpha}})\leq 0.
\end{align}
As a consequence,
\begin{align}
\label{ESTTTTL}
(\|\tilde{u}\|_{X_{\tau,\alpha}}+K\|\tilde{b}\|_{X_{\tau,\alpha}})\langle t\rangle^{1/4-\eta_1}&+\int_0^t\frac{\eta_2}{\langle s\rangle^{1/4+\eta_1}}(\|\tilde{u}(s)\|_{D_{\tau,\alpha}}+K\|\tilde{b}(s)\|_{D_{\tau,\alpha}})ds\nonumber\\
&\qquad\leq (\|\tilde{u}(0)\|_{X_{\tau,\alpha}}+K\|\tilde{b}(0)\|_{X_{\tau,\alpha}})\leq C(1+K)\varepsilon,
\end{align}
where we have used (\ref{IIIIE}).
Then, by noting  that $K=\frac{4C}{\delta}$, one has
\begin{align}
\label{T1}
&\frac{3C_0}{2}(K+1)\int_0^t\langle s\rangle^{-1/4}(\|\tilde{u}(s)\|_{X_{\tau,\alpha}}+\|\tilde{b}(s)\|_{X_{\tau,\alpha}})ds\nonumber\\
=&\frac{3C_0}{2}(\frac{4C}{\delta}+1)\int_0^t\langle s\rangle^{-1/4}(\|\tilde{u}(s)\|_{X_{\tau,\alpha}}+\|\tilde{b}(s)\|_{X_{\tau,\alpha}})ds\nonumber\\
\leq &\frac{3CC_0\varepsilon}{2}(\frac{4C}{\delta}+1)^2\int_0^t\langle s\rangle^{-1/2+\eta_1}ds
\leq 3CC_0\varepsilon(\frac{4C}{\delta}+1)^2\langle t\rangle^{1/2+\eta_1},
\end{align}
and
\begin{align}
\label{T2}
&\frac{3C_0}{2}(K+1)\int_0^t\langle s\rangle^{1/4}(\|\tilde{u}(s)\|_{D_{\tau,\alpha}}+\|\tilde{b}(s)\|_{D_{\tau,\alpha}})ds\nonumber\\
=&\frac{3C_0}{2}(\frac{4C}{\delta}+1)\int_0^t\langle s\rangle^{1/4}(\|\tilde{u}(s)\|_{D_{\tau,\alpha}}+\|\tilde{b}(s)\|_{D_{\tau,\alpha}})ds\nonumber\\
=&\frac{3C_0}{2}(\frac{4C}{\delta}+1)\frac{8}{\delta}\int_0^t\langle s\rangle^{1/2+\eta_1}\frac{\eta_2}{\langle s\rangle^{1/4+\eta_1}}(\|\tilde{u}(s)\|_{D_{\tau,\alpha}}+K\|\tilde{b}(s)\|_{D_{\tau,\alpha}})ds\nonumber\\
\leq&(\frac{4C}{\delta}+1)^2\frac{12CC_0}{\delta}\langle t\rangle^{1/2+\eta_1}\varepsilon.
\end{align}
On the other hand,  (\ref{ODE}) implies that
\begin{align}
\label{ODES}
\tau(t)^{3/2}
=\tau(0)^{3/2}-\frac{3C_0(K+1)}{2}\int_0^t(\langle s\rangle^{-1/4}(\|\tilde{u}\|_{X_{\tau,\alpha}}+\|\tilde{b}\|_{X_{\tau,\alpha}})+\langle s\rangle^{1/4}(\|\tilde{u}\|_{D_{\tau,\alpha}}+\|\tilde{b}\|_{D_{\tau,\alpha}}))ds.
\end{align}
From (\ref{T1}), (\ref{T2}) and (\ref{ODES}), one has
\begin{align*}
\tau(t)^{3/2}\geq \tau_0^{3/2}-\max\{3CC_0(\frac{4C}{\delta}+1)^2\langle t\rangle^{1/2+\eta_1}\varepsilon,\quad (\frac{4C}{\delta}+1)^2\frac{12CC_0}{\delta}\langle t\rangle^{1/2+\eta_1}\varepsilon \},
\end{align*}
for all $t\geq 0$.

Choose $\delta=\frac{1}{\ln(1/\varepsilon)}$. It is straightforward  to show that
\begin{align*}
\tau(t)\geq \frac{\tau_0}{4},
\end{align*}
in the time interval $[0, T_\varepsilon]$, where $T_\varepsilon$ satisfies
\begin{align}
\label{LS}
T_\varepsilon= \bar{C}\left(\frac{1}{\varepsilon(\ln(1/\varepsilon))^3}\right)^{2-4/(\ln(1/\varepsilon)+2)}-1.
\end{align}
This gives the estimate on the lifespan of solution stated  in (\ref{THM3}).

\section{The Proof of Uniqueness Part in Theorem \ref{THM}}
Assume there are two solutions $(\tilde{u}_1, \tilde{b}_1)$ and  $(\tilde{u}_2, \tilde{b}_2)$ to (\ref{3.4}) with the same initial data $(\tilde{u}_0, \tilde{b}_0)$, which satisfies $\|(\tilde{u}_0, \tilde{b}_0)\|_{X_{2\tau_0}, \alpha}\leq \varepsilon$. And the tangential radii of analyticity of $(\tilde{u}_1, \tilde{b}_1)$ and  $(\tilde{u}_2, \tilde{b}_2)$ are $\tau_1(t)$ and $\tau_2(t)$, respectively.

Define $\tau(t)$ by
\begin{align}
\label{U1}
\frac{d(\tau(t))^{3/2}}{dt}+\frac{3C_0(K+1)}{2}&\left(\langle t\rangle^{-1/4}(\|\tilde{u}_1\|_{X_{\tau_1(t),\alpha}}+\|\tilde{b}_1\|_{X_{\tau_1(t),\alpha}})\right.\nonumber\\
&\left.+\langle t\rangle^{1/4}(\|\tilde{u}_1\|_{D_{\tau_1(t),\alpha}}+\|\tilde{b}_1\|_{D_{\tau_1(t),\alpha}})\right)=0,
\end{align}
with initial data
\begin{align}
\label{U2}
\tau(0)=\frac{\tau_0}{8}.
\end{align}
By the estimates obtained  in Section 2, there exists a time interval $[0, T_0]$ with $T_0\leq T_\varepsilon$  such that
\begin{align}
\label{U3}
\frac{\tau_0}{16}\leq \tau(t)\leq \frac{\tau_0}{8}\leq \frac{\min\{\tau_1,\tau_2\}}{2}
\end{align}
for all $t\in[0,T_0]$.

Set  $U=\tilde{u}_1-\tilde{u}_2$ and $B=\tilde{b}_1-\tilde{b}_2$. Then
\begin{align}
\label{U4}
\partial_tU-\partial_y^2U+(u_s+u_1)\partial_xU&+(v_1-v_2)\partial_y\tilde{u}_1-(1+b_1)\partial_xB-(g_1-g_2)\partial_y\tilde{b}_1\nonumber\\
&\qquad\qquad-2\partial_y^2u_s B+(v_1-v_2)\partial_y^2u_s\psi_1+R_{s1}=0,
\end{align}
and
\begin{align}
\label{U5}
\partial_tB-\partial_y^2B-(1+b_1)\partial_xU-(g_1-g_2)\partial_y\tilde{u}_1&+(u_s+u_1)\partial_xB+(v_1-v_2)\partial_y\tilde{b}_1\nonumber\\
&\qquad-(g_1-g_2)\partial_y^2u_s\psi_1+R_{s2}=0,
\end{align}
with the source terms $R_{s1}$ and $R_{s2}$ given by
\begin{align}
\label{U6}
R_{s1}=(u_1-u_2)\partial_x\tilde{u}_2+v_2\partial_yU-(b_1-b_2)\partial_x\tilde{b}_2-g_2\partial_yB+v_2\partial_y^2u_s(\psi_1-\psi_2),
\end{align}
and
\begin{align}
\label{U7}
R_{s2}=-(b_1-b_2)\partial_x\tilde{u}_2-g_2\partial_yU+(u_1-u_2)\partial_x\tilde{b}_2+v_2\partial_yB-g_2\partial_y^2u_s(\psi_1-\psi_2).
\end{align}
Note that the initial data and the boundary conditions are
\begin{align}
\label{DEI}
U(t,x,y)|_{t=0}=0,\qquad B(t, x,y)|_{t=0}=0,
\end{align}
and
\begin{align}
\label{DE}
\left\{
\begin{array}{ll}
U|_{y=0}=0,\\
U|_{y=\infty}=0,
\end{array}
\right.
\qquad\hbox{and}\qquad
\left\{
\begin{array}{ll}
\partial_yB|_{y=0}=0,\\
B|_{y=\infty}=0.
\end{array}
\right.
\end{align}
Similar to Section 2, we have
\begin{align}
\label{U8}
&\frac{d}{dt}\|U\|_{X_{\tau,\alpha}}+\sum_{m\geq 0}\tau^mM_m\frac{\|\theta_\alpha\partial_y\partial_x^m U\|^2_{L^2}}{\|\theta_\alpha\partial_x^m U\|_{L^2}}-\frac{\alpha}{2\langle t\rangle}\|U\|_{X_{\tau,\alpha}}-\frac{C}{\langle t\rangle}\|B\|_{X_{\tau,\alpha}}\\
\leq&\dot{\tau}(t)\|U\|_{Y_{\tau,\alpha}}+\frac{C_0}{(\tau(t))^{1/2}}\left(\langle t\rangle^{-1/4}(\|\tilde{u}_1\|_{X_{\tau(t),\alpha}}+\|\tilde{b}_1\|_{X_{\tau(t),\alpha}})+\langle t\rangle^{1/4}(\|\tilde{u}_1\|_{D_{\tau(t),\alpha}}+\|\tilde{b}\|_{D_{\tau_1(t),\alpha}})\right)\nonumber\\
&\times (\|U\|_{Y_{\tau,\alpha}}+\|B\|_{Y_{\tau,\alpha}})\nonumber\\
&+\frac{C_0}{(\tau(t))^{1/2}}(\|\tilde{u}_2\|_{Y_{\tau,\alpha}}+\|\tilde{b}_2\|_{Y_{\tau,\alpha}})
(\langle t\rangle^{-1/4}(\|U\|_{X_{\tau,\alpha}}+
\|B\|_{X_{\tau,\alpha}})+\langle t\rangle^{1/4}(\|U\|_{D_{\tau,\alpha}}+\|B\|_{D_{\tau,\alpha}}))\nonumber,
\end{align}
and
\begin{align}
\label{U9}
&\frac{d}{dt}\|B\|_{X_{\tau,\alpha}}+\sum_{m\geq 0}\tau^mM_m\frac{\|\theta_\alpha\partial_y\partial_x^m U\|^2_{L^2}}{\|\theta_\alpha\partial_x^m U\|_{L^2}}-\frac{\alpha}{2\langle t\rangle}\|B\|_{X_{\tau,\alpha}}\\
\leq&\dot{\tau}(t)\|B\|_{Y_{\tau,\alpha}}+\frac{C_0}{(\tau(t))^{1/2}}\left(\langle t\rangle^{-1/4}(\|\tilde{u}_1\|_{X_{\tau(t),\alpha}}+\|\tilde{b}_1\|_{X_{\tau(t),\alpha}})+\langle t\rangle^{1/4}(\|\tilde{u}_1\|_{D_{\tau(t),\alpha}}+\|\tilde{b}_1\|_{D_{\tau(t),\alpha}})\right)\nonumber\\
&\times (\|U\|_{Y_{\tau,\alpha}}+\|B\|_{Y_{\tau,\alpha}})\nonumber\\
&+\frac{C_0}{(\tau(t))^{1/2}}(\|\tilde{u}_2\|_{Y_{\tau,\alpha}}+\|\tilde{b}_2\|_{Y_{\tau,\alpha}})
(\langle t\rangle^{-1/4}(\|U\|_{X_{\tau,\alpha}}+
\|B\|_{X_{\tau,\alpha}})+\langle t\rangle^{1/4}(\|U\|_{D_{\tau,\alpha}}+\|B\|_{D_{\tau,\alpha}}))\nonumber.
\end{align}
Then, we have
\begin{align}
\label{U10}
&\frac{d}{dt}(\|U\|_{X_{\tau,\alpha}}+K\|B\|_{X_{\tau,\alpha}})+\sum_{m\geq 0}\tau^mM_m\left(\frac{\|\theta_\alpha\partial_y\partial_x^m U\|^2_{L^2}}{\|\theta_\alpha\partial_x^m U\|_{L^2}}+K\frac{\|\theta_\alpha\partial_y\partial_x^m B\|^2_{L^2}}{\|\theta_\alpha\partial_x^m B|_{L^2}}\right)\nonumber\\
&-\frac{\alpha}{2\langle t\rangle}\|U\|_{X_{\tau,\alpha}}-\frac{2C+K\alpha}{2\langle t\rangle}\|B\|_{X_{\tau,\alpha}}\nonumber\\
\leq&\left(\dot{\tau}(t)+\frac{C_0(1+K)}{(\tau(t))^{1/2}}\left(\langle t\rangle^{-1/4}(\|\tilde{u}_1\|_{X_{\tau(t),\alpha}}+\|\tilde{b}_1\|_{X_{\tau(t),\alpha}})+\langle t\rangle^{1/4}(\|\tilde{u}_1\|_{D_{\tau(t),\alpha}}+\|\tilde{b}_1\|_{D_{\tau(t),\alpha}})\right)\right.\nonumber\\
&\left.\times(\|U\|_{Y_{\tau,\alpha}}+K\|B\|_{Y_{\tau,\alpha}})\right)\\
&+\frac{C_0(1+K)}{(\tau(t))^{1/2}}(\|\tilde{u}_2\|_{Y_{\tau,\alpha}}+\|\tilde{b}_2\|_{Y_{\tau,\alpha}})
(\langle t\rangle^{-1/4}(\|U\|_{X_{\tau,\alpha}}+
\|B\|_{X_{\tau,\alpha}})+\langle t\rangle^{1/4}(\|U\|_{D_{\tau,\alpha}}+\|B\|_{D_{\tau,\alpha}}))\nonumber.
\end{align}
From (\ref{U1}), one has
\begin{align}
\label{U11}
\dot{\tau}(t)+\frac{C_0(K+1)}{(\tau(t))^{1/2}}&\left(\langle t\rangle^{-1/4}(\|\tilde{u}_1\|_{X_{\tau,\alpha}}+
\|\tilde{b}_1\|_{X_{\tau,\alpha}})\right.\nonumber\\
&\left.+\langle t\rangle^{1/4}(\|\tilde{u}_1\|_{D_{\tau,\alpha}}+\|\tilde{b}_1\|_{D_{\tau,\alpha}})\right)(\|U\|_{Y_{\tau,\alpha}}+\|B\|_{Y_{\tau,\alpha}})
\leq 0,
\end{align}
because $\tau(t)\leq \tau_1(t)$ and the norms
$X_{\tau, \alpha}$ and $D_{\tau, \alpha}$  are increasing in $\tau$.

By the inequalities (\ref{ESTU}), (\ref{ESTB}) and (\ref{U11}), one has
\begin{align}
\label{U12}
&\frac{d}{dt}(\|U\|_{X_{\tau,\alpha}}+K\|B\|_{X_{\tau,\alpha}})+\frac{\alpha(1-2\beta_1)}{2\langle t\rangle}\|U\|_{X_{\tau,\alpha}}
+\frac{\alpha(1-2\beta_2)-\frac{2C}{K}}{2\langle t\rangle}K\|B\|_{X_{\tau,\alpha}}\nonumber\\
&+\frac{\alpha^{1/2}\beta_1}{2\langle t\rangle^{1/2}}\|U\|_{D_{\tau,\alpha}}+\frac{\alpha^{1/2}\beta_2}{2\langle t\rangle^{1/2}}K\|B\|_{D_{\tau,\alpha}}\\
\leq &\frac{C_0(1+K)}{(\tau(t))^{1/2}}(\|\tilde{u}_2\|_{Y_{\tau,\alpha}}+\|\tilde{b}_2\|_{Y_{\tau,\alpha}})
(\langle t\rangle^{-1/4}(\|U\|_{X_{\tau,\alpha}}+
\|B\|_{X_{\tau,\alpha}})+\langle t\rangle^{1/4}(\|U\|_{D_{\tau,\alpha}}+\|B\|_{D_{\tau,\alpha}}))\nonumber,
\end{align}
for $\beta_1, \beta_2\in (0, 1/2)$.
Since
\begin{align}
\label{U13}
\|\tilde{u}_2\|_{Y_{\tau,\alpha}}\leq \frac{1}{\tau}\|\tilde{u}_2\|_{X_{2\tau,\alpha}}\leq \frac{1}{\tau}\|\tilde{u}_2\|_{X_{\tau_2,\alpha}}
\leq \frac{C(1+K)}{\tau}\varepsilon\langle t\rangle^{-1/4+\eta_1}
\end{align}
and
\begin{align}
\label{U14}
\|\tilde{b}_2\|_{Y_{\tau,\alpha}}\leq \frac{1}{\tau}\|\tilde{b}_2\|_{X_{2\tau,\alpha}}\leq \frac{1}{\tau}\|\tilde{b}_2\|_{X_{\tau_2,\alpha}}
\leq \frac{C(1+K)}{\tau}\varepsilon\langle t\rangle^{-1/4+\eta_1},
\end{align}
we have
\begin{align}
\label{U15}
\frac{C_0(1+K)}{(\tau(t))^{1/2}}(\|\tilde{u}_2\|_{Y_{\tau,\alpha}}+\|\tilde{b}_2\|_{Y_{\tau,\alpha}})\leq \frac{2(1+K)^2CC_0\varepsilon}{(\tau(t))^{3/2}\langle t\rangle^{1/4-\eta_1}}.
\end{align}
Notice that $t\in [0,T_\varepsilon]$ with $T_\varepsilon=\varepsilon^{-2+\delta_0}$, where $\delta_0$ is a fixed
small positive constant. As in Section 3, we can choose $\alpha=1/2-\delta, \beta_1=\beta_2=\frac{\delta}{2}, K=\frac{4C}{\delta}$ and
$\delta=1/\ln(1/\varepsilon)$, then $\eta_1$ can be chosen to be $\delta$.  Let $\varepsilon$ suitably small to have
\begin{align*}
\frac{\alpha(1-2\beta_1)}{2}>\frac{2(1+K)^2CC_0\varepsilon\langle t\rangle^{1/2+\eta_1}}{(\tau(t))^{3/2}},\qquad
\frac{\alpha(1-2\beta_1)-2C/K}{2}>\frac{2(1+K)^2CC_0\varepsilon\langle t\rangle^{1/2+\eta_1}}{(\tau(t))^{3/2}K},
 \end{align*}
and
\begin{align*}
\frac{\alpha^{1/2}\beta_1}{2}>\frac{2(1+K)^2CC_0\varepsilon\langle t\rangle^{1/2+\eta_1}}{(\tau(t))^{3/2}},\qquad
\frac{\alpha^{1/2}\beta_2}{2}>\frac{2(1+K)^2CC_0\varepsilon\langle t\rangle^{1/2+\eta_1}}{(\tau(t))^{3/2}K}.
  \end{align*}
 (\ref{U12}) and (\ref{U15}) imply that
\begin{align}
\label{U16}
\frac{d}{dt}(\|U\|_{X_{\tau,\alpha}}+K\|B\|_{X_{\tau,\alpha}})+\eta_3(\|U\|_{X_{\tau,\alpha}}
+K\|B\|_{X_{\tau,\alpha}})\leq 0
\end{align}
for suitably small $\eta_3>0$ and any $t\in [0,T_\varepsilon]$.  It implies uniqueness of solution to (\ref{3.4})  in the time interval $[0, T_\varepsilon]$.\\

\noindent{\bf Acknowledgement:} Feng Xie' research  was supported by National Nature Science Foundation of
China 11571231, the China Scholarship
Council and Shanghai Jiao Tong University SMC(A). Tong Yang's research was supported by internal research funding of City University of Hong Kong, 7004847.

\end{document}